\documentclass[12pt]{article}
\usepackage{amsbsy,amsfonts,amsmath,amssymb,amsthm,bbm,euscript,eepic,rotating}
\usepackage{graphicx}
\def\UseSection{
        \numberwithin{equation}{section}
    \theoremstyle{plain}
        \newtheorem{theorem}    {Theorem}[section]
        \DefineTheorems 
}

\def\DefineTheorems{
    
    \newtheorem{lemma}      [theorem] {Lemma}
    
    \newtheorem{prop}       [theorem] {Proposition}
    
    \newtheorem{cor}        [theorem] {Corollary}

    \theoremstyle{definition}
    \newtheorem{defn}       [theorem] {Definition}

    \theoremstyle{definition}

}

\newcommand{\bt}   {\begin{theorem}}
\newcommand{\et}   {\end  {theorem}}
\newcommand{\bl}   {\begin{lemma}}
\newcommand{\el}   {\end  {lemma}}
\newcommand{\bp}   {\begin{prop}}
\newcommand{\ep}   {\end  {prop}}
\newcommand{\bc}   {\begin{cor}}
\newcommand{\ec}   {\end  {cor}}
\newcommand{\bd}   {\begin{defn}}
\newcommand{\ed}   {\end  {defn}}

\newcommand{\ba}   {\begin{array}}
\newcommand{\ea}   {\end  {array}}
\newcommand{\be}   {\begin{enumerate}}
\newcommand{\ee}   {\end  {enumerate}}
\newcommand{\bi}   {\begin{itemize}}
\newcommand{\ei}   {\end  {itemize}}

\def\eq#1\en{\begin{equation}#1\end{equation}}
\def\eqsplit#1\ensplit{
    \begin{equation}\begin{split}#1\end{split}\end{equation}
    }
\def\eqalign#1\enalign{
    \begin{align}#1\end{align}
    }
\def\eqmul#1\enmul{
    \begin{multline}#1\end{multline}
    }
\newcommand{\eqarrstar} {\begin{eqnarray*}}
\newcommand{\enarrstar} {\end{eqnarray*}}
\newcommand{\eqarray}   {\begin{eqnarray}}
\newcommand{\enarray}   {\end{eqnarray}}

\newcommand{\lbeq}[1]  {\label{e:#1}}
\newcommand{\refeq}[1] {\eqref{e:#1}}    

\makeatletter
\newcommand{\labelcounter}[2]{{%
    \stepcounter{#1}
    \protected@write\@auxout{}%
    {\string\newlabel{#2}{{\csname the#1\endcsname}{\thepage}}}%
    {\ref{#2}}
    }}
\makeatother
\newcommand{\sss}   { \scriptscriptstyle }

\newcommand{\Rbold} {{\mathbb R}}

\newcommand{\Zbold} {{\mathbb Z}}

\newcommand{\ovec}  {\boldsymbol{o}}

\newcommand{\uvec}  {\boldsymbol{u}}
\newcommand{\vvec}  {\boldsymbol{v}}

\newcommand{\xvec}  {\boldsymbol{x}}
\newcommand{\yvec}  {\boldsymbol{y}}
\newcommand{\zvec}  {\boldsymbol{z}}

\newcommand{\Scal}   {\mathcal{S}}

\newcommand{\Rd}    {{ {\Rbold}^d}}
\newcommand{\Zd}    {{ {\Zbold}^d }}

\newcommand{\spose}[1] {{\hbox to 0pt{#1\hss}} }
\newcommand{\ltapprox} {\mathrel{\spose{\lower 3pt\hbox{$\mathchar"218$}}
 \raise 2.0pt\hbox{$\mathchar"13C$}}}
\newcommand{\gtapprox} {\mathrel{\spose{\lower 3pt\hbox{$\mathchar"218$}}
 \raise 2.0pt\hbox{$\mathchar"13E$}}}

\UseSection  
\theoremstyle{plain}
\newtheorem{thm}{Theorem}[section]

\newcommand{\bb}{\underline{b}}
\newcommand{\bC}{\tilde C}

\newcommand{\daw}{\downarrow}
\newcommand{\db}{\Longrightarrow}

\newcommand{\BDcup}[1]{~\underset{#1}{\Dot{\bigcup}}~}
\newcommand{\dc}{d_{\rm c}}

\newcommand{\Elow}{E_{\prec}}
\newcommand{\Ehigh}{E_{\succ}}
\newcommand{\ind}[1]{\mathbbm{1}{\scriptstyle\{#1\}}}
\newcommand{\indic}[1]{\ind{#1}}

\newcommand{\mE}{{\mathbb E}}

\newcommand{\mP}{{\mathbb P}}

\newcommand{\mR}{{\mathbb R}}
\newcommand{\mZ}{{\mathbb Z}}

\newcommand{\nnmb}{\nonumber}

\newcommand{\Rp}{\mR_+}

\newcommand{\scscst}{\scriptscriptstyle}

\newcommand{\tb}{\overline{b}}
\newcommand{\tbp}
 {\tb^{\raisebox{-2pt}{\scriptsize$\prime$}}}

\newcommand{\too}{\longrightarrow}

\newcommand{\chicp}{\chi^{\scscst\rm cp}}
\newcommand{\chiop}{\chi^{\scscst\rm op}}
\newcommand{\chiperc}{\chi^{\scscst\rm pe}}
\newcommand{\chisaw}{\chi^{\scscst\rm sa}}
\newcommand{\pc}{p_{\rm c}}
\newcommand{\pccp}{\pc^{\scscst\rm cp}}
\newcommand{\pcop}{\pc^{\scscst\rm op}}
\newcommand{\pcperc}{\pc^{\scscst\rm pe}}
\newcommand{\pcsaw}{\pc^{\scscst\rm sa}}
\newcommand{\Pisaw}{\Pi^{\scscst\rm sa}}
\newcommand{\hPisaw}{\hat\Pi^{\scscst\rm sa}}
\newcommand{\taucp}{\tau^{\scscst\rm cp}}
\newcommand{\tauop}{\tau^{\scscst\rm op}}
\newcommand{\tauperc}{\tau^{\scscst\rm pe}}
\newcommand{\tausaw}{\tau^{\scscst\rm sa}}

\newcommand{\uaw}{\uparrow}

\newcommand{\vep}{\varepsilon}

\newcommand{\zd}{\mZ^d}
\newcommand{\Zp}{\mZ_+}

\newcommand{\udb}{\Longleftrightarrow}
\newcommand{\cntd}{\longleftrightarrow}
\newcommand{\nn}{\nonumber}

 \newcommand{\smallsup}[1] {{\scriptscriptstyle{({#1}})}}
 
 \newcommand{\R}{\Rbold}
\newcommand{\Z}{\Zbold}

\newcommand{\ctx}[1] { \xleftrightarrow{#1} }
\newcommand{\xleftrightarrow}[2][]{\leftarrow\hspace{-2.0ex}\xrightarrow[#1]{#2}}

\newcommand{\sN}{{\sss N}}
\newcommand{\sR}{{\sss R}}

\title{Critical points for spread-out self-avoiding walk, percolation\\
and the contact process above the upper critical dimensions}

\author{
    Remco van der Hofstad\footnote{Department of Mathematics and
        Computer Science, Eindhoven University of Technology,
        5600 MB Eindhoven, The Netherlands.
        {\tt rhofstad@win.tue.nl}}\\
    Akira Sakai\footnote{EURANDOM, P.O. Box 513, 5600 MB Eindhoven,
        The Netherlands.  {\tt sakai@eurandom.tue.nl}}
    }

\date{December 1, 2003\footnote{Revised August 3, 2004}}

\begin{document}
\setlength{\baselineskip}{16pt}
\maketitle


\begin{abstract}
We consider self-avoiding walk and percolation in $\Zd$, oriented percolation
in $\Zd\times\Zp$, and the contact process in $\Zd$, with $p\,D(\,\cdot\,)$
being the coupling function whose range is denoted by $L<\infty$.  For
percolation, for example, each bond $\{x,y\}$ is occupied with probability
$p\,D(y-x)$.  The above models are known to exhibit a phase transition when
the parameter $p$ varies around a model-dependent critical point $\pc$.  We
investigate the value of $\pc$ when $d>6$ for percolation and $d>4$ for the
other models, and $L\gg1$.  We prove in a unified way that
$\pc=1+C(D)+O(L^{-2d})$, where the universal term 1 is the mean-field
critical value, and the model-dependent term $C(D)=O(L^{-d})$ is written
explicitly in terms of the random walk transition probability $D$.  We also
use this result to prove that $\pc=1+cL^{-d}+O(L^{-d-1})$, where $c$ is a
model-dependent constant.  Our proof is based on the lace expansion for each
of these models.
\end{abstract}

\section{Introduction and main results}\label{s:results}
Self-avoiding walk, percolation, and the contact process are well-known
models that exhibit critical phenomena.  For percolation in two or higher
dimensions, for example, there exists a percolation threshold $\pcperc$
such that there is almost surely no infinite cluster for $p<\pcperc$,
while for $p>\pcperc$ there is almost surely a unique infinite cluster.
As $p\uaw\pcperc$, the average cluster size and the correlation length
diverge.  The precise value of $\pcperc$ depends on the details of the
model, and is only explicitly known in a few cases, such as for
two-dimensional nearest-neighbor bond percolation \cite{k80}.

In this paper, we will consider self-avoiding walk, percolation, oriented
percolation and the contact process, where the interaction range $L$ is
taken to be large.  When $L\gg1$, the interaction in the considered models
is relatively weak, and therefore the critical values can be expected to be
close to the critical value 1 of the respective mean-field models, i.e.,
random walk and branching random walk.  We study the difference of the
critical values and 1 for the above four models as $L\to\infty$.  It turns
out that, above the respective upper critical dimensions, we can write this
difference to leading order as simple functions of the underlying random
walk.

\subsection{Models}\label{ss:models}
First, we define the models.  A self-avoiding walk is a path $\omega$ in
the $d$-dimensional integer lattice $\Zd$ with $\omega(i)\ne\omega(j)$
for every distinct $i,j\in\{0,1,\dots,|\omega|\}$.  We also take the 
zero-step walk into account.  We define the weight of a non-zero path 
$\omega$ by
\begin{align}\lbeq{weight-def}
W_p(\omega)=p^{|\omega|}\prod_{i=1}^{|\omega|}D\big(\omega(i)-\omega
 (i-1)\big),
\end{align}
where $D$ is a probability distribution on $\Zd$, and let $W_p(\omega)=1$ 
if $|\omega|=0$.  We suppose that $D$ is symmetric with respect to the 
lattice symmetries and that $D(o)=0,$ where $o$ is the origin in $\Z^d$.  
A more detailed definition will be given below.
The self-avoiding walk two-point function is defined by
\begin{align}\lbeq{saw-two}
\tausaw_p(x)=\sum_{\substack{\omega:o\too x\\ {\rm saw}}}W_p(\omega),
\end{align}
where the sum is over all self-avoiding paths from $o$ to $x$.  It is known
(see, e.g., \cite{ms93}) that there is a critical value $\pcsaw$
such that
\begin{align}\lbeq{chisaw}
\chisaw_p=\sum_{x\in\Zd}\tausaw_p(x)
\end{align}
is finite if and only if $p<\pcsaw$ and diverges as $p\uaw\pcsaw$.

For percolation, each bond $\{x,y\}$ is occupied with probability
$p\,D(y-x)$ and vacant with probability $1-p\,D(y-x)$, independently of
the other bonds, where $p\in[0,\|D\|_\infty^{-1}]$.  Since $\sum_xD(x)=1$,
the percolation parameter $p$ is the expected number of occupied bonds per
site.  We denote by $\mP_p$ the probability distribution of the bond
variables.  We say that $x$ is connected to $y$, and write $x\cntd y$,
if either $x=y$ or there is a path of occupied bonds between $x$ and $y$.
The percolation two-point function and its sum over $\Zd$ are denoted by
\begin{align}\lbeq{perc-two}
\tauperc_p(x)=\mP_p(o\longleftrightarrow x),&&
\chiperc_p=\sum_{x\in\Zd}\tauperc_p(x).
\end{align}
Similarly to self-avoiding walk, there is a critical value $\pcperc$ such
that $\chiperc_p$ is finite if and only if $p<\pcperc$ and diverges as
$p\uaw\pcperc$ (see, e.g., \cite{g99}).

Oriented percolation is a time-directed version of percolation.  Each bond
$((x,t),(y,t+1))$ is an ordered pair of sites in $\Zd\times\Zp$,
and is occupied with probability $p\,D(y-x)$ and vacant with probability
$1-p\,D(y-x)$, independently of the other bonds, where
$p\in[0,\|D\|_\infty^{-1}]$.  We say that $(x,s)$ is connected to $(y,t)$,
and write $(x,s)\too(y,t)$, if either $(x,s)=(y,t)$ or there is an oriented
path of occupied bonds from $(x,s)$ to $(y,t)$.  Let $\mP_p$ be the
probability distribution of the bond variables.  The oriented percolation
two-point function and its sum over $\Zd\times\Zp$ are denoted by
\begin{align}\lbeq{op-two}
\tauop_p(x,t)=\mP_p((o,0)\too (x,t)),&&
\chiop_p=\sum_{t\in\Zp}\,\sum_{x\in\Zd}\tauop_p(x,t).
\end{align}
Also oriented percolation exhibits a phase transition such that
$\chiop_p<\infty$ if and only if $p$ is less than the critical value $\pcop$,
and that $\chiop_p\uaw\infty$ as $p\uaw\pcop$ (see, e.g., \cite{gh01}).

The contact process is a model of the spread of an infection in $\Zd$,
and is a continuous-time version of oriented percolation in $\Zd\times\Rp$. We
now describe a graphical representation for the contact process.  Along each time line
$\{x\}\times\Rp$, where $x\in \Z^d$, we place points according to a Poisson process with
intensity 1, independently of the other time lines.  For each ordered pair
of distinct time lines from $\{x\}\times\Rp$ to $\{y\}\times\Rp$, we place
oriented bonds $((x,t),(y,t))$, $t\geq0$, according to a Poisson process with
intensity $p\,D(y-x)$, independently of the other Poisson processes, where
the parameter $p\geq0$ is the infection rate.  We say that $(x,s)$ is connected to $(y,t)$,
and write $(x,s)\too(y,t)$, if either $(x,s)=(y,t)$ or there is an oriented
path in $\Zd\times\Rp$ from $(x,s)$ to $(y,t)$ using the Poisson bonds and
time-line segments traversed in the increasing-time direction without
traversing the Poisson points.  Let $\mP_p$ be the corresponding probability
distribution.  We denote the contact process two-point function and its
integro-sum over $\Zd\times\Rp$ by
\begin{align}\lbeq{cp-two}
\taucp_p(x,t)=\mP_p((o,0)\too(x,t)),&&
\chicp_p=\int_0^\infty\!\!dt~\sum_{x\in\Zd}\taucp_p(x,t).
\end{align}
Again there is a critical value $\pccp$ such that
$\chicp_p$ is finite if and only if $p<\pccp$ and diverges as $p\uaw\pccp$
(see, e.g., \cite{Ligg99}).

We will omit the superscript referring to the precise model, and write $\pc$
when referring to the critical values in all models simultaneously.  The goal
in this paper is to study $\pc$ when the range $L$ of $D$ is sufficiently
large.  In the proofs, we will have versions of $D$ in mind which are such
that $L^d D(Lx)$ is a discrete approximation of a function on $\R^d$.  We
will formalize this assumption on $D$ in the following definition:

\begin{defn}\label{defn:D}
Let $h$ be a probability distribution over $\Rd\setminus\{o\}$, which is
invariant under rotations by $\pi/2$ and reflections in the coordination
hyperplanes.  We suppose that $h$ is piecewise continuous, so that
$\int_{\mR^d}h(x)\,d^dx\equiv1$ can be approximated by the Riemann sum
$L^{-d}\sum_{x\in\Zd}h(x/L)$.  Then, we define
\begin{align}
D(x)=\frac{h(x/L)}{\sum_{y\in\Zd}h(y/L)}.
\end{align}
\end{defn}

We will make heavy use of results proved elsewhere for the models under
consideration.  For these results, some further assumptions are made on $D$,
of which we now list the most important ones.  We require that there exist
$c>0$, $C<\infty$, $\eta\in(0,1)$ such that
\begin{align}\lbeq{Dprop}
\sup_{x\in\Zd}D(x)\leq CL^{-d},&&
\eta\wedge(cL^2|k|^2)\leq1-\hat D(k)\leq2-\eta,
\end{align}
where $\hat D(k)=\sum_{x\in\Zd}D(x)\,e^{ik\cdot x}$ and
$|k|^2=\sum_{j=1}^dk_j^2$.  There are a few more minor requirements
that depend on the precise model under investigation.  For details, see
\cite{HS90a} for percolation and \cite{hsa04,hs02,hs01,hs03} for the other
three models, for which the requirements are virtually identical.  A simple
example of $D$, where all the above assumptions are satisfied, is
\begin{align}\lbeq{dex}
D(x)=\frac{\indic{0<\|x\|_\infty\leq L}}{(2L+1)^d-1},
\end{align}
for which $h(x)=2^{-d}$ if $0<\|x\|_\infty\leq1$ and $h(x)=0$ otherwise.

We denote by $D*G$ the convolution of $D$ and a function $G$ in $\mZ^d$,
and by $D^{*n}$ the $n$-fold convolution of $D$ in $\Zd$, where we define
$D^{*0}(x)=\delta_{o,x}$.  We will frequently use
    \begin{align}\lbeq{d-gauss}
    D^{*n}(x)\leq\delta_{0,n}\delta_{o,x}+\frac{O(\beta)}{(1\vee n)^{d/2}},
    \end{align}
where
    \begin{align}
    \beta=L^{-d}.\lbeq{beta-def}
    \end{align}
The inequality \refeq{d-gauss} is a consequence of \refeq{Dprop},
as we will show in Appendix~\ref{s:lclt}.

\subsection{Main results}
Let $\dc$ denote the respective upper critical dimensions, i.e., $\dc=6$ for
percolation and $\dc=4$ for the other three models.  In this paper, we
investigate the respective critical values when $d>\dc$ and $L\gg1$, in a
unified fashion.

\begin{thm}\label{thm:main}
For each model with $d>\dc$, as $L\rightarrow \infty$,
\begin{align}
\pcsaw,~\pccp
&=1+\sum_{n=2}^\infty D^{*n}(o)+O(\beta^2),\lbeq{sawcp-critpt}\\
\pcop &=1+\frac12\sum_{n=2}^\infty D^{*2n}(o)+O(\beta^2),
 \lbeq{op-critpt}\\
\pcperc &=1+D^{*2}(o)+\frac12\sum_{n=3}^\infty(n+1)\,D^{*n}(o)
 +O(\beta^2).\lbeq{perc-critpt}
\end{align}
\end{thm}

The universal term 1 is the critical value for the mean-field models
(random walk and branching random walk).  Note that, by \refeq{d-gauss},
the model-dependent terms in \refeq{sawcp-critpt}--\refeq{perc-critpt}
are $O(\beta)$.  In Section~\ref{ss:overview}, we will intuitively
explain why the model-dependent terms have the above respective forms.


We next compute the dependence on $\beta$ more explicitly, and compute the
coefficients of $\beta$ in $\pc-1$ explicitly.  For this, we let $U$ be the
uniform probability distribution over $[-1,1]^d\subset\Rd$, i.e., for
$x\in\Rd$,
\begin{align}\lbeq{U-def}
U(x)=2^{-d}\indic{\|x\|_\infty\leq1},
\end{align}
and denote by $U^{\star n}$ the $n$-fold convolution of $U$ in $\Rd$.
Then, the leading order coefficient in $\beta$ for $\pc$ is given in
the following theorem:

\begin{thm}\label{thm:uniform}
Fix $D$ as in \refeq{dex}, and let $d>\dc$.  As $L\rightarrow\infty$,
\begin{align}
\pcsaw,~\pccp&=1+\beta\sum_{n=2}^\infty U^{\star n}(o)+O(\beta L^{-1}),
 \lbeq{sawcp-uniform}\\
\pcop &=1+\frac\beta2\sum_{n=2}^\infty U^{\star2n}(o)+O(\beta L^{-1}),\\
\pcperc &=1+\beta\bigg[U^{\star2}(o)+\frac12\sum_{n=3}^\infty(n+1)\,
 U^{\star n}(o)\bigg]+O(\beta L^{-1}).
\end{align}
\end{thm}

We now comment on the relation between the asymptotics in
Theorems~\ref{thm:main}--\ref{thm:uniform}.  The advantage of
Theorem~\ref{thm:uniform} is that it is more concrete, and the continuum
limit of the critical points appears explicitly.  However, the error term
in Theorem~\ref{thm:main} is $O(\beta^2)$, while in Theorem~\ref{thm:uniform}
it is equal to $O(\beta L^{-1})$, which is much larger.  In order to compute
the critical value more precisely, Theorem~\ref{thm:main} gives a much more
powerful result, at the expense of having to compute the random walk terms
appearing in its statement.  In principle, it should be possible to compute
the coefficients of $\beta L^{-1}, \beta L^{-2}, \dots, \beta L^{-d+1}$, but
this requires a substantial amount of work.  Finally, it should be possible
to compute the random walk sums in Theorem~\ref{thm:main} for other examples
than the one in \refeq{dex}, but we refrain from doing so.

We now summarize previous results on the critical values.  We start with
self-avoiding walk.  Penrose's result in \cite{p94} implies that the
critical value for self-avoiding walk defined by
\refeq{dex} with $L\gg1$ satisfies
    \begin{align}\lbeq{p94}
    1+c\,\beta^{2/7}\log\beta^{-1}\geq\pcsaw\geq\begin{cases}
    1,&\text{if }~d\geq3,\\
    1+c'\beta\log\beta^{-1},&\text{if }~d=2,\\
    1+c''\beta^{4/5},&\text{if }~d=1,
    \end{cases}
    \end{align}
for some $\beta$-independent constants $c,c',c''$.  For {\it spread-out
lattice trees}, a related result with a different leading term,
namely $=e^{-1}$, was also obtained in \cite{p94}.  For $d>4$, Madras and Slade
\cite[Corollary 6.2.7]{ms93} improved \refeq{p94} to $\pcsaw=1+O(\beta)$.
In \cite{hs02,hs03}, this result was extended to more general $D$ as
defined in Definition~\ref{defn:D}.  We will rely on the results in
\cite{hs02,hs03}, whose proof is based on the {\it lace expansion} and a
generalized inductive approach.  We will also use the lace expansion to
derive the expression of the $O(\beta)$ term in \refeq{sawcp-critpt}.

For percolation, the best previous result is
$\pcperc=1+O(\beta^{2/d-\epsilon})$ for $d>6$ and $L\gg1$, where
$\epsilon>0$ is an arbitrarily small number \cite{HHS01a}.  However,
if we combine Lemma~\ref{lem-T} proved below and the estimates for
the lace expansion in \cite{HS90a}, then we obtain the better estimate
$\pcperc=1+O(\beta)$.  The result in \refeq{perc-critpt}, which is also
obtained by an application of the lace expansion, identifies the
expression of this $O(\beta)$ term.

When $d>4$ and $L\gg1$, both $\pcop$ and $\pccp$ were proved to be
$1+O(\beta)$ \cite{hsa04,hs02,hs01}.  Similarly to self-avoiding walk,
the proofs of these results rely on the lace expansion and an adaptation
of the inductive approach.  The contact process defined in terms of $D$
of \refeq{dex} was first considered by Bramson, Durrett and Swindle
\cite{bds89}, and they proved that, as $L\to\infty$,
    \begin{align}\lbeq{bds}
    \pccp-1\asymp f(\beta)\equiv\begin{cases}
    \beta,&\text{if }d\geq3,\\
    \beta\log\beta^{-1},&\text{if }d=2,\\
    \beta^{2/3},&\text{if }d=1,
    \end{cases}
    \end{align}
where $\pccp-1\asymp f(\beta)$ means that the ratio $(\pccp-1)/f(\beta)$ is
bounded away from zero and infinity.  Later, Durrett and Perkins \cite{dp99}
proved that
    \begin{align}\lbeq{ed}
    \lim_{L\to\infty}\frac{\pccp-1}{f(\beta)}=\begin{cases}
    \sum_{n=2}^\infty U^{\star n}(o),&\text{if }d\geq3,\\
    3/(2\pi),&\text{if }d=2.
    \end{cases}
    \end{align}
Our result \refeq{sawcp-uniform} in Theorem~\ref{thm:uniform} is stronger
when $d>4$ in the sense that not only the coefficient of $\beta$, but also
the speed of convergence in \refeq{ed} is identified.  In \cite{hsa04}, we
also obtained certain lace expansion results for a local mean-field limit,
where the range and time grow large simultaneously, for the contact process
in $d\leq4$, and we expect that these results could be used to prove
a stronger version of \refeq{ed} for $d=3,4$, as well as for oriented
percolation.  However, this will need serious work using block constructions
as used in \cite{dp99}.

We expect that \refeq{sawcp-critpt}--\refeq{perc-critpt} remain valid for
$d=\dc-1$ and $\dc$ when we change $O(\beta^2)$ to $o(\beta)$.  As mentioned
above, this is the case for the contact process \cite{dp99}.  When
$d\leq\dc-2$, the second terms in \refeq{sawcp-critpt}--\refeq{perc-critpt}
diverge, so that Theorem~\ref{thm:main} cannot hold.  However, we
expect that the asymptotics of the critical point will, as for the contact
process, again be described by the divergence of the sums in
\refeq{sawcp-critpt}--\refeq{perc-critpt}.

When $d>\dc$, we expect that the $O(\beta^2)$ terms could be identified in
terms of $D$ as well, using a similar method as in this paper, but to
do so will require a serious amount of work.

A related problem is to obtain the asymptotics of the critical points for the
nearest-neighbor models, when $D(x)=(2d)^{-1}\indic{|x|=1}$ and $d\to\infty$.
In \cite{hs95}, $\pcsaw$ was proved to have an asymptotic expansion into
powers of $(2d)^{-1}$, and the first six coefficients were obtained.  For
unoriented percolation, the first three coefficients were computed in
\cite{hs95} and \cite{HS04a}, but the proof of the asymptotic expansion only
appeared in \cite{HS04b}.  The proofs of these results are again based on the lace
expansion.  For nearest-neighbour oriented percolation and the
nearest-neighbour contact process, it is proved that $\pcop=1+O(d^{-2})$
(see \cite{cd83}) and $\pccp=1+O(d^{-1})$ (see, e.g., \cite{l97}), using
different methods.

\subsection{Overview of the proof}\label{ss:overview}
To prove Theorem~\ref{thm:main}, we will apply the lace expansion (see, e.g.,
\cite{HS90a,hsa04,hs03,ms93,ny93}).  For example, the lace expansion for
self-avoiding walk gives the recurrence relation
    \begin{align}\lbeq{saw-renewal}
    \tausaw_p(x)=\delta_{o,x}+\sum_{v}[p\,D(v)+\Pisaw_p(v)]\,
    \tausaw_p(x-v),
    \end{align}
where $\Pisaw_p(x)$ is a certain expansion coefficient.  It was proved
in \cite{hs02,hs03} that $\hPisaw_p\equiv\sum_x\Pisaw_p(x)=O(\beta)$ for
$p\leq\pcsaw$, if $d>4$ and $L\gg1$ (see Section~\ref{s:saw}).  Summing
both sides of \refeq{saw-renewal} over $x\in\Zd$ and solving the resulting
equation in terms of $\chisaw_p$, we obtain
    \begin{align}\lbeq{saw-chi}
    \chisaw_p=(1-p-\hPisaw_p)^{-1},
    \end{align}
and thus
    \begin{align}\lbeq{saw-pc}
    \pcsaw=1-\hPisaw_{\pcsaw}.
    \end{align}
To estimate $\pcsaw$, we thus need to investigate $\hPisaw_{\pcsaw}$.
We will prove that, since $\pcsaw=1+O(\beta)$, we can replace
$\hPisaw_{\pcsaw}$ by $\hPisaw_1$ up to an error of order $O(\beta^2)$.
When $p=1$, the exponentially growing factor $p^{|\omega|}$ in
\refeq{weight-def} does not play any role, and $\hPisaw_1$ can be
investigated in terms of random walks.  This is the key ingredient for
the proof of Theorem~\ref{thm:main}.


The strategy for percolation models is the same as above.  There is
a similar recursion relation to \refeq{saw-renewal}, with some
model-dependent expansion coefficient $\Pi_p(x)$.  Therefore, to obtain
the formulae in Theorem~\ref{thm:main}, we will have to investigate
$\hat\Pi_1=\sum_x\Pi_1(x)$, again in terms of random walks.

As we will explain in Sections~\ref{s:saw}--\ref{s:perc}, $\hPisaw_1$
and $\hat\Pi_1$ can be described by an alternating sum of a model-dependent
sequence $\hat\pi_1^{\sss(N)}$ for $N\geq0$, where $\hat\pi_1^{\sss(N)}$
for $N\geq1$ decays as $\beta^N$ for all models.  For self-avoiding walk,
$\hat\pi_1^{\sss(0)}$ equals zero, while $\hat\pi_1^{\sss(0)}$ for
percolation models is nearly a half of $\hat\pi_1^{\sss(1)}$.  (This is why
we have the factor $\frac12$ in \refeq{op-critpt}--\refeq{perc-critpt}.)
Therefore,
roughly speaking, we only need to investigate $\hat\pi_1^{\sss(1)}$ to
obtain \refeq{sawcp-critpt}--\refeq{perc-critpt}.  We will show later that
the diagrammatic interpretation of $\hat\pi_1^{\sss(1)}$ for self-avoiding
walk is a single random walk taking more than one step and going back
to the starting point (cf., \refeq{sawcp-critpt}), while the diagrammatic
interpretation of $\hat\pi_1^{\sss(1)}$ for percolation models is that two
random walks, at least one of which is non-vanishing, meet at some point.
Therefore, the correction to the mean-field value 1 are related to
{\it random walk loops}.

For loops in the time-oriented percolation models, the lengths in
the time-increasing direction of these two walks have to be equal
(which explains the sum over even convolution powers in \refeq{op-critpt}),
while for unoriented percolation this is not the case (which explain the
sum over all powers and the factor $n+1$ in \refeq{perc-critpt}).

For the contact process, the two paths are {\it continuous time} random walk paths,
for which the number of convolution powers of $D$ is equal to the number of
spatial steps made by the random walk up to a given time, which has
a Poisson distribution. Therefore, the sum over the
convolution powers of $D$ is not restricted to even powers, and we see that
the correction to the mean-field value for the contact process and
oriented percolation are different. For the contact process, it will turn
out that also the factor $\frac12$ in \refeq{op-critpt} disappears, which
is due to the fact that the two walks are in fact {\it avoiding} each other,
and which will be explained in more detail in Section \ref{ss:op}. This is an
intuitive explanation of the model-dependent terms in
\refeq{sawcp-critpt}--\refeq{perc-critpt}.

We organize the rest of this paper as follows.  We begin with self-avoiding
walk in Section~\ref{s:saw}, and explain the key steps to estimate
$\pcsaw$.  Following the same steps, we estimate $\pcop$ and $\pccp$ in
Section~\ref{ss:op}, and $\pcperc$ in Section~\ref{ss:perc}.  Finally,
we prove an extension of \refeq{d-gauss} in Appendix~\ref{s:lclt}, and
Theorem~\ref{thm:uniform} in Appendix~\ref{s:uniform}.

\section{Critical point for self-avoiding walk}\label{s:saw}
In this section, we prove \refeq{sawcp-critpt}, using \refeq{saw-pc}.
Throughout this section, we will omit the superscript ``sa'' and write,
e.g., $\pc=\pcsaw$ and $\hat\Pi_p=\hPisaw_p$.

Before computing the asymptotics of $\hat\Pi_{\pc}$ in \refeq{saw-pc},
we first note that $\pc\geq1$.  This is because the removal of the self-avoidance
constraint in \refeq{saw-two} results in $\sum_{\omega:o\too x}W_p(\omega)$,
whose sum over $x\in\Zd$ equals $(1-p)^{-1}$ for any $p\leq1$.
For self-avoiding walk,
    \begin{align}\lbeq{saw-lace}
    \Pi_p(x)=\sum_{N=1}^\infty(-1)^N\pi_p^{\scscst(N)}(x),
    \end{align}
where, e.g., $\pi_p^{\scscst(1)}(x)$ is a ``1-loop diagram'' at the origin
\cite{ms93}:
    \begin{align}\lbeq{saw-1stdiag}
    \pi_p^{\scscst(1)}(x)=\delta_{o,x}~(pD*\tau_p)(o)=\delta_{o,x}\sum_{
    \substack{\omega:o\too o\\ |\omega|\geq1}}W_p(\omega)\,I(\omega),
    \end{align}
where $I(\omega)=1$ if there are no self-intersection points except
for $\omega(0)=\omega(|\omega|)$, otherwise $I(\omega)=0$.

For $d>4$ and $L\gg1$, it was proved in \cite{hs03} that, for
$\hat\pi_p^{\scscst(N)}=\sum_x\pi_p^{\scscst(N)}(x)$, we have
    \begin{align}\lbeq{saw-piNbd}
    \hat\pi_p^{\scscst(N)}\leq O(\beta)^N,&& p\,\partial_p\hat\Pi_p\leq O(\beta),
    \end{align}
for all $p\leq\pc$ and $N\geq1$.  Together with \refeq{saw-pc} and
\refeq{saw-lace}, we immediately obtain that $\pc=1+O(\beta)$.  Moreover, by the
mean-value theorem, there is a $p\in(1,\pc)$ such that
    \begin{align}\lbeq{saw-pcrewrite}
    \pc=1-\hat\Pi_1-(\hat\Pi_{\pc}-\hat\Pi_1)&=1-\hat\Pi_1-(\pc-1)~\partial_p
    \hat\Pi_p=1+\hat\pi_1^{\scscst(1)}+O(\beta^2),
    \end{align}
where
    \begin{align}
    \hat\pi_1^{\scscst(1)}=\sum_{\substack{\omega:o\too o\\ |\omega|\geq1}}
    W_1(\omega)\,I(\omega)
    =\sum_{n=2}^\infty D^{*n}(o)-\sum_{\substack{\omega:o\too o\\ |\omega|
    \geq1}}W_1(\omega)~[1-I(\omega)].\lbeq{saw-pi1exp}
    \end{align}

To complete the proof of \refeq{sawcp-critpt}, it thus suffices to prove
that the second term in the right-hand side of \refeq{saw-pi1exp} is
$O(\beta^2)$ if $d>4$.  We first note that $I(\omega)$ is an indicator function.
If $I(\omega)=0$, so that $1-I(\omega)=1$, then there must be a pair
$\{s,t\}\ne\{0,|\omega|\}$ with $0\leq s<t\leq|\omega|$ such that
$\omega(s)=\omega(t)$.  Denoting the parts of $\omega$ corresponding to
these three time intervals by $\omega_i$, $i=1,2,3$, respectively, we obtain
    \begin{align}\lbeq{saw-error}
    \sum_{\substack{\omega:o\too o\\ |\omega|\geq1}}W_1(\omega)~[1-I(\omega)]
    \leq\sum_{x\in\Zd}\,\sum_{\substack{\omega_1,\omega_3:o\too x\\ |\omega_1|
    +|\omega_3|\geq1}}\,\sum_{\substack{\omega_2:x\too x\\ |\omega_2|\geq1}}\,
    \prod_{i=1}^3W_1(\omega_i)=(D*G^{*2})(o)~(D^{*2}\!*G)(o),
    \end{align}
where $G(x)=\sum_{n=0}^\infty D^{*n}(x)$, and
$(D^{*2}*G)(o)=\sum_{n=2}^\infty D^{*n}(o)=O(\beta)$ if $d>2$.  Moreover,
by \refeq{d-gauss},
    \eq
    (D*G^{*2})(o)=\sum_{n=1}^\infty nD^{*n}(o)=O(\beta)
    \en
if $d>4$.  This completes the proof of \refeq{sawcp-critpt}
for self-avoiding walk.
\qed

\section{Critical points for percolation models}\label{s:perc}
In this section, we compute the asymptotics of the critical values
for the other three models, and thus complete the proof of
Theorem~\ref{thm:main}.

To discuss oriented percolation and the contact process simultaneously,
it is convenient to introduce the following oriented percolation on
$\Zd\times\vep\Zp$, which is the time-discretized contact process with
a discretization parameter $\vep\in(0,1]$.  A bond is a directed pair
$((x,t),(y,t+\vep))$ of sites in $\Zd\times\vep\Zp$.  Each bond
is either occupied or vacant, independently of the other bonds, and a bond
$((x,t),(y,t+\vep))$ is occupied with probability
\begin{align}\lbeq{bprob}
q_p(y-x)=\begin{cases}
 1-\vep,&\text{if }x=y,\\
 p\vep D(y-x),&\text{if }x\ne y,
 \end{cases}
\end{align}
provided that $\sup_x q_p(x)\leq1$.  In this notation, the model with
$\vep=1$ is the usual oriented percolation model as defined in
Section \ref{ss:models}, and the weak limit as
$\vep\daw0$ is the contact process \cite{bg91}.  Similarly to
oriented percolation with $\vep=1$, for each $\vep \in (0,1]$,
there is a critical value $\pc^{\sss(\vep)}$
for every $\vep\in(0,1],$ such that $\pc^{\sss(1)}=\pcop$ and
$\lim_{\vep\downarrow 0}\pc^{\sss(\vep)}=\pccp$ \cite{s01}. We will call the
model with $\vep\in (0,1]$ the {\it time-discretized contact process}.

To summarise notation for percolation and the time-discretized contact
process, we will write $\Lambda=\Zd$ for percolation and $\Lambda=
\Zd\times\vep\Zp$ for oriented percolation. For notational convenience,
we will take $\vep=1$ for percolation. We will also use bold letters
to represent elements of $\Lambda$.  For example, $\ovec=o$, $\xvec=x$
for percolation, and $\ovec=(o,0)$, $\xvec=(x,t)$ for the time-discretized
contact process. For a bond $b=(\uvec,\vvec)$, we write $\bb=\uvec$ and
$\tb=\vvec$.  We also omit the superscripts $\vep$, pe, op and cp, if no
confusion can arise.

As mentioned in Section~\ref{s:results}, the lace expansion for percolation
models takes a similar form as in \refeq{saw-renewal}, and reads (see, e.g.,
\cite{HS90a, hsa04})
    \begin{align}\lbeq{perc-renewal}
    \tau_p(\xvec)=[\delta_{\ovec,\xvec}+\Pi_p(\xvec)]+\sum_{\uvec,\vvec\in\Lambda}
    [\delta_{\ovec,\uvec}+\Pi_p(\uvec)]\,q_p(\vvec-\uvec)\,\tau_p(\xvec-\vvec).
    \end{align}
In particular, $q_p(\vvec-\uvec)=p\,D(v-u)$ for percolation and oriented
percolation for which $\vep=1$.  (To unify notation, we recall that we regard
unoriented percolation as a model with $\vep=1$.)  The lace expansion
coefficient $\Pi_p(\xvec)$ equals
    \begin{align}\lbeq{perc-lace}
    \Pi_p(\xvec)=\sum_{N=0}^\infty(-1)^N\pi_p^{\scscst(N)}(\xvec),
    \end{align}
where $\pi_p^{\scscst(N)}(\xvec)$, $N\geq0$, are model-dependent diagram
functions.  The result of the lace expansion will be explained in
Sections~\ref{ss:op}--\ref{ss:perc}.  For the time-discretized contact process with
$\vep\in(0,1]$, $d>\dc$ and $L\gg1$, it has been proved \cite{hsa04,hs01}
that $\hat\Pi_p\equiv\vep\sum_{\xvec\in\Lambda}\Pi_p(\xvec)$ is
$O(\beta)\,\vep^2$ for all $p\leq\pc$.  The same estimate is proved to hold
for unoriented percolation (with $\vep=1$), using the lace expansion in
\cite{HS90a} and Lemma~\ref{lem-T} proved below in Section~\ref{ss:perc}.

As in the derivation of \refeq{saw-chi}, solving \refeq{perc-renewal} in
terms of $\chi_p=\vep\sum_{\xvec\in\Lambda}\tau_p(\xvec)$ gives
    \begin{align}
    \chi_p=\frac{1+\frac1\vep\hat\Pi_p}{1-p-(1-\vep+p\vep)\frac1{\vep^2}
    \hat\Pi_p},
    \end{align}
and thus, equating the denominator to zero,
    \begin{align}\lbeq{perc-pc}
    \pc=1-\frac1{\vep^2}\hat\Pi_{\pc}-(\pc-1)\frac1\vep\hat\Pi_{\pc}.
    \end{align}
This expression holds uniformly in $\vep$. We will use it to compute $\pcop$
and $\pcperc$ by taking $\vep=1$ and $\pccp$ by taking the limit when
$\vep\daw0$ \cite{s01}, respectively.  In particular, the third term is
$O(\beta^2)$ when $\vep=1$, and it has no contribution in the limit
$\vep\daw0$. Therefore, we are left to prove that, apart from an error term
of order $O(\beta^2)$, the second term in \refeq{perc-pc} equals the second
term in \refeq{sawcp-critpt} when $\vep\daw0$, and equals the second term in
\refeq{op-critpt} for oriented percolation and that in \refeq{perc-critpt}
for (unoriented) percolation when $\vep=1$.  We again note that
$\pc^{\smallsup{\vep}}\geq1$, since
$\chi_p\leq\vep\sum_{n=0}^\infty\sum_xq_p^{*n}(x)=(1-p)^{-1}$ for $p\leq1$.
In addition, similarly to \refeq{d-gauss}, when $p=1$ and $\vep<1$, we have
    \begin{align}\lbeq{q-gauss}
    q^{*n}(x)\equiv q_1^{*n}(x)\leq(1-\vep)^n\,\delta_{o,x}+\frac{O(\beta)}
    {[1\vee(n\vep)]^{d/2}}.
    \end{align}
We will prove \refeq{q-gauss} in Appendix~\ref{s:lclt}.  Note that when $\vep=1$,
\refeq{q-gauss} reduces to \refeq{d-gauss}.

To complete the proof of Theorem~\ref{thm:main}, we investigate
$\hat\Pi_{\pc}$ for oriented percolation and the contact process
in Section~\ref{ss:op}, and for unoriented percolation in
Section~\ref{ss:perc}.

\subsection{Asymptotics of $\pcop$ and $\pccp$}
\label{ss:op}
In this section, we investigate $\hat\Pi_{\pc}$ for the discretized contact
process, and derive \refeq{op-critpt} for oriented percolation (i.e.,
$\vep=1$) and \refeq{sawcp-critpt} for the contact process (i.e., $\vep\daw0$).

To describe the diagram functions $\pi_p^{\sss(N)}(\xvec)$, $N\geq0$, we
need some definitions.  We say that $\xvec$ is {\it doubly connected to}
$\yvec$, if either $\xvec=\yvec$ or there are at least two nonzero
bond-disjoint occupied paths from $\xvec$ to $\yvec$.  Following the
notation in \cite{hs01} as closely as possible, we denote this event by
$\xvec\db\yvec$, and define
    \begin{align}\lbeq{op-pi0}
    \hat \pi_p^{\scscst(0)}=\vep \sum_{\xvec \in \Lambda} \pi_p^{\scscst(0)}(\xvec),
    \qquad \text{where}\qquad \pi_p^{\scscst(0)}(\xvec)=\mP_p(\ovec\db\xvec)-\delta_{\ovec,\xvec}.
    \end{align}
If $\ovec$ is connected but not doubly connected to $\xvec$, there is
a {\it pivotal bond} $b=(\bb,\tb)$ for $\ovec\too\xvec$ such that both
$\ovec\too\bb$ and $\tb\too\xvec$ occur, and that $\ovec\too\xvec$ occurs
if and only if $b$ is set occupied.  For $A\subseteq\Lambda$, we say that
$\yvec$ {\it is connected to $\xvec$ through} $A$ when every occupied path
from $\yvec\too\xvec$ has at least one bond with an endpoint in $A$.  We define
$E(b,\xvec;A)$ to be the event that $b$ is occupied, that $\tb\too\xvec$
through $A$, and that there are no pivotal bonds $b'$ for $\tb\too\xvec$
such that $\tb\too\bb'$ through $A$.  Let $\bC^b(\ovec)$ be the set of
vertices in $\Lambda$ connected from $\ovec$ without using $b$.  Then,
    \begin{align}
    \hat \pi_p^{\scscst(1)}=\vep \sum_{\xvec \in \Lambda} \pi_p^{\scscst(1)}(\xvec),
    \qquad \text{where}\qquad\pi_p^{\scscst(1)}(\xvec)=\sum_b\mP_p\big(\ovec\db\bb;\,E(b,\xvec;\bC^b
    (\ovec))\big).\lbeq{op-pi1}
    \end{align}
The higher order diagram functions $\pi_p^{\scscst(N)}(\xvec)$,
$N\geq2$, are defined in a similar way, but are irrelevant in this
paper (see \cite[Section~3]{hsa04} for a complete definition,
with slightly different notation).

For $d>4$ and $L\gg1$, it was proved in \cite{hsa04} that, for
$\hat\pi_p^{\scscst(N)}=\vep\sum_{\xvec\in\Lambda}\pi_p^{\scscst(N)}(\xvec)$,
we have
    \begin{align}\lbeq{op-piNbd}
    \hat\pi_p^{\scscst(N)}\leq O(\beta)^{N\vee1}\,\vep^2,&&
    p\,\partial_p\hat\Pi_p\leq O(\beta)\,\vep^2,
    \end{align}
for all $p\leq\pc$ and $N\geq0$.  Together with \refeq{perc-lace} and
\refeq{perc-pc}, we obtain $\pc=1+O(\beta)$.  Moreover, by the mean-value
theorem, there is a $p\in(1,\pc)$ such that
    \begin{align}\lbeq{op-pcrewrite}
    \pc=1-\frac1{\vep^2}\hat\Pi_{\pc}-(\pc-1)\frac1\vep\hat\Pi_{\pc}
    &=1-\frac1{\vep^2}\hat\Pi_1-(\pc-1)\frac1{\vep^2}\partial_p\hat\Pi_p
    +O(\beta^2)\,\vep\nnmb\\
    &=1-\frac1{\vep^2}\hat\pi_1^{\scscst(0)}+\frac1{\vep^2}\hat\pi_1^{\scscst(1)}
    +O(\beta^2).
    \end{align}
To prove \refeq{sawcp-critpt}--\refeq{op-critpt}, it thus suffices to
investigate $\hat\pi_1^{\scscst(0)}$ and $\hat\pi_1^{\scscst(1)}$.

\begin{proof}[Analysis of $\hat\pi_1^{\scscst(0)}$]
We prove
    \begin{align}\lbeq{pi0-main}
    \frac1{\vep^2}\hat\pi_1^{\scscst(0)}\begin{cases}
    =\frac12\sum_{n=2}^\infty D^{*2n}(o)+O(\beta^2),&\text{for }~\vep=1,\\
    \to\sum_{n=2}^\infty D^{*n}(o)+O(\beta^2),&\text{when }~\vep\daw0.
    \end{cases}
    \end{align}

Recall \refeq{op-pi0}.  To describe a double connection by a pair of two
random walk paths, we order the support of $D$ in an arbitrary but fixed
manner.  For $x,y$ in the support of $D$, we write $x\prec y$ if $x$ is
{\it lower} than $y$ in that order.  For a pair of paths consisting of
bonds in $\Lambda$, $\omega=(b_1,\dots,b_{\sN})$ and $\omega'=(b'_1,\dots,b'_{\sN})$
with $\bb_1=\bb'_1$ and $\tb_{\sN}=\tbp_{\sN}$, we say that $\omega$ is {\it lower}
than $\omega'$, denoted by $\omega\prec\omega'$, if at the first time
$n\in\{1,\dots,N\}$ when $\omega$ is incompatible with $\omega'$
(therefore $b_i=b'_i$ for all $i<n$) we have $\tb_n\prec\tbp_n$.
We also say that $\omega_2$ is {\it higher} than $\omega_1$.

A path $\omega=(b_1,\dots,b_{|\omega|})$ is said to be {\it occupied} if
all bonds along $\omega$ are occupied.  We define $\Elow(\omega)$ to be
the event that $\omega$ is the lowest occupied path among all occupied
paths from $\bb_1$ to $\tb_{|\omega|}$, and that there is another occupied
path $\omega'$ from $\bb_1$ to $\tb_{|\omega|}$ which is bond-disjoint
from $\omega$ (denoted by $\omega\cap\omega'=\varnothing$).  Given a path
$\omega$, we also define $\Ehigh(\omega';\omega)$ to be the event that
$\omega'$ is the highest occupied path among all occupied paths from $\bb_1$
to $\tb_{|\omega|}$ that are bond-disjoint from $\omega$.  Such an occupied
path $\omega'$ exists on $\{\bb_1\db\tb_{|\omega|}\}\cap\Elow(\omega)$ by
definition.

Using the above notation, we have, for $\xvec\ne\ovec$,
    \begin{align}\lbeq{pi0-eventreexpr}
    \{\ovec\db\xvec\}=\BDcup{\substack{\omega_1,\omega_2:\ovec\too\xvec\\
    \omega_1\cap\omega_2=\varnothing\\ \omega_1\prec\omega_2}}\big\{
    \omega_1,\omega_2\text{ occupied};\,\Elow(\omega_1)\cap\Ehigh
    (\omega_2;\omega_1)\big\}.
    \end{align}
We define the right-hand side to be empty if $\xvec=\ovec$.  Then,
    \begin{align}\lbeq{pi-reexpr}
    \hat\pi_1^{\scscst(0)}=\vep\sum_{\xvec\in\Lambda}~\sum_{\substack{\omega_1,
    \omega_2:\ovec\too\xvec\\ \omega_1\cap\omega_2=\varnothing\\ \omega_1\prec
    \omega_2}}\mP_1\big(\omega_1,\omega_2\text{ occupied};\,\Elow(\omega_1)\cap
    \Ehigh(\omega_2;\omega_1)\big).
    \end{align}
Since $\mP_1$ is a product measure, if we ignore
$\Elow(\omega_1)\cap\Ehigh(\omega_2;\omega_1)$, then we obtain
    \begin{align}
    &\sum_{\substack{\omega_1,\omega_2:\ovec\too\xvec\\ \omega_1\cap\omega_2=
    \varnothing\\ \omega_1\prec\omega_2}}\mP_1(\omega_1,\omega_2\text{ occupied})
    \nnmb\\
    &\quad=\sum_{\substack{u,v:u\prec v\\ y,z:y\ne z}}q(u)\,q(v)\,q(x-y)\,q(x-z)
    \sum_{\substack{\omega_1:\uvec\too\yvec\\ \omega_2:\vvec\too\zvec\\ \omega_1
    \cap\omega_2=\varnothing}}\mP_1(\omega_1\text{ occupied})~\mP_1(\omega_2
    \text{ occupied}),\lbeq{pi-premain}
    \end{align}
where $\uvec=(u,\vep)$, $\vvec=(v,\vep)$, $\yvec=(y,t-\vep)$,
$\zvec=(z,t-\vep)$, and $q(x)=q_1(x)$ (cf., \refeq{q-gauss}).  By an
inclusion-exclusion relation, 
the correction is bounded by
    \[ \sum_{\substack{\omega_1,\omega_2:\ovec\too\xvec\\ \omega_1\cap\omega_2=
    \varnothing\\ \omega_1\prec\omega_2}}\big[\mP_1\big(\omega_1,\omega_2
    \text{ occupied};\,\Elow(\omega_1)^{\rm c}\big)+\mP_1\big(\omega_1,\omega_2
    \text{ occupied};\,\Ehigh(\omega_2;\omega_1)^{\rm c}\big)\big]. \]
We will prove below that, for $E$ equal to $\Elow(\omega_1)$ or
$\Ehigh(\omega_2;\omega_1)$,
    \begin{align}\lbeq{intersect1}
    \vep\sum_{\xvec\in\Lambda}~\sum_{\substack{\omega_1,\omega_2:\ovec\too\xvec\\
    \omega_1\cap\omega_2=\varnothing}}\mP_1(\omega_1,\omega_2\text{ occupied};\,
    E^{\rm c})=O(\beta^2)\,\vep^2.
    \end{align}

We investigate \refeq{pi-premain} to obtain the expression of
$O(\beta)$ from \refeq{pi-reexpr}.  If we ignore the restriction
$\omega_1\cap\omega_2=\varnothing$, then we obtain
    \begin{align}\lbeq{pi-main}
    \sum_{\substack{u,v:u\prec v\\ y,z:y\ne z}}q(u)\,q(v)\,q(x-y)\,
    q(x-z)~q^{*(t/\vep-2)}(y-u)~q^{*(t/\vep-2)}(z-v),
    \end{align}
where $t/\vep\in[2,\infty)\cap\Zp$.  We will prove below that the correction
satisfies
    \begin{align}\lbeq{intersect2}
    \vep\sum_{\xvec\in\Lambda}\,\sum_{\substack{u,v:u\prec v\\ y,z:y\ne z}}\!
    \!q(u)\,q(v)\,q(x-y)\,q(x-z)\!\!\sum_{\substack{\omega_1:\uvec\too\yvec\\
    \omega_2:\vvec\too\zvec\\ \omega_1\cap\omega_2\ne\varnothing}}\!\!\mP_1(
    \omega_1\text{ occupied})~\mP_1(\omega_2\text{ occupied})=O(\beta^2)\,\vep^2.
    \end{align}
Therefore, we only need to consider the contribution to \refeq{pi-reexpr}
from \refeq{pi-main}.
By changing variables as $y'=x-y$ and $z'=x-z$ and using the symmetry between
$u\prec v$ and $u\succ v$, the sum of \refeq{pi-main} over $x\in\Zd$ equals
    \begin{gather}
    \sum_{\substack{u,v:u\prec v\\ y',z':y'\ne z'}}q(u)\,q(v)\,q(y')\,q(z')\sum_x
    q^{*(t/\vep-2)}(x-y'-u)~q^{*(t/\vep-2)}(x-z'-v)\nnmb\\
    =\frac12\sum_{\substack{u,v:u\ne v\\ y,z:y\ne z}}q(u)\,q(v)\,q(y)\,q(z)~
    q^{*(2t/\vep-4)}(v+z-y-u).\lbeq{decomp}
    \end{gather}
Recall \refeq{bprob}.  Since there is at most one temporal (or vertical)
bond growing out of every site in $\Lambda$, we must have
$q(u)=\vep D(u)$ or $q(v)=\vep D(v)$, so that we obtain at
least one factor of $\vep$.  By the same reason, we should have
$q(y)=\vep D(y)$ or $q(z)=\vep D(z)$, so that we obtain
a second factor of $\vep$.  Therefore, the number of
combinations for the product of four factors of $q$ in \refeq{decomp} is
nine: one combination is proportional to $\vep^4$, four others are
proportional to $(1-\vep)\,\vep^3$, and the remaining four are proportional
to $(1-\vep)^2\,\vep^2$.  Only the first case arises for oriented percolation
for which $\vep=1$ , while only the third case arises for the contact process
for which $\vep\daw0$, respectively.

We first complete the proof of \refeq{pi0-main} for oriented percolation.
When $\vep=1$, and using inclusion-exclusion on the restrictions
$u\neq v$ and $y\neq z$, the sum of \refeq{decomp} over $t\geq2$ equals
    \begin{gather}\lbeq{oppi-main}
    \frac12\sum_{u,v,y,z}D(u)\,D(v)\,D(y)\,D(z)\sum_{t=2}^\infty D^{*(2t-4)}
    (v+z-y-u)+O(\beta^2)=\frac12\sum_{t=2}^\infty D^{*2t}(o)+O(\beta^2),
    \end{gather}
where we use \refeq{Dprop} to obtain an error of order $O(\beta^2)$ that
comes from contributions where $u=v$ or $y=z$.

For the contact process, for which $\vep\daw0$, the leading contribution is
due to the four combinations of order $(1-\vep)^2\,\vep^2$ mentioned above,
where either $u$ or $v$ is $o$, and either $y$ or $z$ is $o$.  Therefore,
the coefficient of $(1-\vep)^2\,\vep^2$ in \refeq{decomp} is
    \begin{align*}
    &\frac12\bigg[\sum_{u,y}D(u)\,D(y)~q^{*(2t/\vep-4)}(-y-u)+\sum_{u,z}D(u)
    \,D(z)~q^{*(2t/\vep-4)}(z-u)\\
    &\quad+\sum_{v,y}D(v)\,D(y)~q^{*(2t/\vep-4)}(v-y)+\sum_{v,z}D(v)\,D(z)
    ~q^{*(2t/\vep-4)}(v+z)\bigg]=2\big(D^{*2}\!*q^{*(2t/\vep-4)}\big)(o).
    \end{align*}
Summing this expression (multiplied by $\vep$) over
$t/\vep\in[2,\infty)\cap\Zp$ gives
    \begin{align}
    2\int_{\Box_\pi}\frac{d^dk}{(2\pi)^d}~\hat D(k)^2~\vep\sum_{n=0}^\infty
    \big[1-\vep+\vep\hat D(k)\big]^{2n}&=\int_{\Box_\pi}\frac{d^dk}{(2
    \pi)^d}~\frac{2\hat D(k)^2}{[1-\hat D(k)][2-\vep+\vep\hat D(k)]}\nnmb\\
    &\stackrel{\vep\daw0}{\too}\int_{\Box_\pi}\frac{d^dk}{(2\pi)^d}~\frac{
    \hat D(k)^2}{1-\hat D(k)}=\sum_{n=2}^\infty D^{*n}(o),\lbeq{cppi-main}
    \end{align}
where $\Box_\pi=[-\pi,\pi]^d$.  This completes the proof of \refeq{pi0-main}.
\end{proof}


\begin{proof}[Analysis of $\hat\pi_1^{\scscst(1)}$]
We prove that $\frac1{\vep^2}\hat\pi_1^{\scscst(1)}$ is asymptotically
twice as large as the right-hand side of \refeq{pi0-main}:
\begin{align}\lbeq{pi1-main}
\frac1{\vep^2}\hat\pi_1^{\scscst(1)}\begin{cases}
 =\sum_{n=2}^\infty D^{*2n}(o)+O(\beta^2),&\text{for }~\vep=1,\\
 \to2\sum_{n=2}^\infty D^{*n}(o)+O(\beta^2),&\text{when }~\vep\daw0.
 \end{cases}
\end{align}

For a bond $b$, let $\{b\db\xvec\}$ be the event that $b$ is occupied and
$\tb\db\xvec$.  We define $\{\uvec\too b\}$ and a joint event
$\{\uvec\too b\db\xvec\}$ similarly.  For events $E_1$ and $E_2$, we denote
by $E_1\circ E_2$ the event that $E_1$ and $E_2$ occur {\it disjointly},
i.e., using disjoint bond sets of bonds (see e.g., \cite[Section 2.3]{g99}).
Recalling \refeq{op-pi1} and distinguishing between $\bb=\ovec$ and $\bb\ne\ovec$, we
can rewrite $\hat\pi_1^{\scscst(1)}$ as
    \begin{gather}
    \hat\pi_1^{\scscst(1)}=\vep\sum_{\uvec,\xvec\in\Lambda}\mP_1\big(\{(\ovec,
    \uvec)\too\xvec\}\circ\{\ovec\too\xvec\}\big)+\vep\sum_{\xvec\in\Lambda}\,
    \sum_{b:\bb\ne\ovec}\mP_1\big(\ovec\db\bb;\,E(b,\xvec;\bC^b(\ovec))\big)
    \nnmb\\
    -\vep\sum_{\uvec,\xvec\in\Lambda}\mP_1\Big(\big\{\{(\ovec,\uvec)\too\xvec\}
    \circ\{\ovec\too\xvec\}\big\}\setminus E((\ovec,\uvec),\xvec;\bC^{\sss
    (\ovec,\uvec)}(\ovec))\Big).\lbeq{pi1-reexpr1}
    \end{gather}
We will extract the leading contribution from the first term.  Note that
$\{(\ovec,\uvec)\too\xvec\}\circ\{\ovec\too\xvec\}$ is almost identical to
$\{\ovec\db\xvec\}=\{\ovec\too\xvec\}\circ\{\ovec\too\xvec\}$.  However, the
symmetry between the two connections from $\ovec$ to $\xvec$ is
lost in the former event, due to the bond $(\ovec,\uvec)$. We will use this
symmetry breaking in a convenient manner. Recall that below \refeq{pi0-main},
the support of $D$ was ordered in an arbitrary way. Now, instead, we choose
the ordering such that, for $\uvec=(u,\vep)$, the element $u$ in the support of
$D$ is {\it minimal}. This will ensure that the lowest occupied
path $\omega_1$ from $\ovec$ to $\xvec$ will use the bond $(\ovec,\uvec)$.
We also write $\Elow^{\uvec}(\omega_1)$ and $\Ehigh^{\uvec}(\omega_2;\omega_1)$
for $\Elow(\omega_1)$ and $\Ehigh(\omega_2;\omega_1)$ in this $\uvec$-dependent
ordering. Therefore, (cf., \refeq{pi0-eventreexpr}),
    \begin{align}\lbeq{pi1-eventreexpr}
    \{(\ovec,\uvec)\too\xvec\}\circ\{\ovec\too\xvec\}=\BDcup{\substack{\omega_1:
    (\ovec,\uvec)\too\xvec\\ \omega_2:\ovec\too\xvec\\ \omega_1\cap\omega_2=
    \varnothing}}&\big\{\omega_1,\omega_2\text{ occupied};\,\Elow^{\uvec}(\omega_1)
    \cap\Ehigh^{\uvec}(\omega_2;\omega_1)\big\},
    \end{align}
and its contribution to \refeq{pi1-reexpr1} is
    \begin{align}\lbeq{pi1-reexpr2}
    \vep\sum_{\xvec\in\Lambda}~\sum_{\substack{\omega_1:(\ovec,\uvec)\too
    \xvec\\ \omega_2:\ovec\too\xvec\\ \omega_1\cap\omega_2=\varnothing}}
    &\mP_1\big(\omega_1,\omega_2\text{ occupied};\,\Elow^{\uvec}(\omega_1)
    \cap\Ehigh^{\uvec}(\omega_2;\omega_1)\big),
    \end{align}
where $\omega_1:(\ovec,\uvec)\too\xvec$ is a path from $\ovec$ to
$\xvec$ starting by the bond $(\ovec,\uvec)$.  Ignoring the condition
$\Elow^{\uvec}(\omega_1)\cap\Ehigh^{\uvec}(\omega_2;\omega_1)$
as in \refeq{pi-premain} and following the same strategy as
in estimating $\hat\pi_1^{\scscst(0)}$,
we obtain the main contribution to \refeq{pi1-main}.  The leading term
of $\frac1{\vep^2}\hat\pi_1^{\scscst(1)}$ is twice as large as that of
$\frac1{\vep^2}\hat\pi_1^{\scscst(0)}$, because the symmetry is broken and
we do not obtain the factor $\frac 12$ as in \refeq{decomp}
(cf., \refeq{pi-reexpr} and \refeq{pi1-reexpr2}).

To complete the proof of \refeq{pi1-main}, it suffices to show that the
second and third terms in \refeq{pi1-reexpr1} are both $O(\beta^2)\,\vep^2$.
The event in the second term of \refeq{pi1-reexpr1} implies the existence of
$\yvec\in\Lambda$ such that $\{\ovec\too\yvec\too\bb\}\circ\{\ovec\too\bb\}$
and $\{\yvec\too\xvec\}\circ\{b\too\xvec\}$ occur disjointly.  Let $\omega_1$
denote a path from $\ovec$ to $\xvec$ through $\yvec$, $\omega_2$ denote
another path from $\ovec$ to $\xvec$ via the the bond $b$ with $\bb=\zvec$,
and $\omega_3$ denote another path from $\yvec$ to $\zvec$.  Then, the second
term in \refeq{pi1-reexpr1} is bounded by
\begin{align}\lbeq{pi1-error}
\vep\sum_{\substack{\xvec,\yvec,\zvec\in\Lambda\\ \zvec\ne\ovec,\xvec}}\,
\sum_{\substack{\omega_1:\ovec\too\yvec\too\xvec\\ \omega_2:\ovec\too\zvec
\too\xvec\\ \omega_3:\yvec\too\zvec\\ \omega_i\cap\omega_j=\varnothing,~i
\ne j}}\,\prod_{i=1}^3\mP_1(\omega_i\text{ occupied}),
\end{align}
since $\mP_1$ is a product measure.  The third term in \refeq{pi1-reexpr1}
is also bounded by the above expression.  This is because the event in the
third term in \refeq{pi1-reexpr1} implies existence of $\yvec\in\Lambda$
and a pivotal bond $b=(\zvec,\cdot\,)$ for $\uvec\too\xvec$ such that
$\{\ovec\too\yvec\too\xvec\}$, $\{(\ovec,\uvec)\too b\too\xvec\}$ and
$\{\yvec\too\zvec\}$ occur disjointly.  We thus obtain \refeq{pi1-error}
by the same random walk representation.

Therefore, it is sufficient to prove that \refeq{pi1-error} is bounded by
$O(\beta^2)\,\vep^2$.  When $\vep=1$, we simply ignore the restriction
$\omega_i\cap\omega_j=\varnothing$, $i\ne j$, and apply the Gaussian bound
\refeq{d-gauss} to the part of $\omega_1$ from $\yvec$ to $\xvec$ and to
the part of $\omega_2$ from $\ovec$ to $\zvec$.  Since $\yvec\ne\xvec$ and
$\zvec\ne \ovec$, the term $\delta_{\ovec,\xvec}$ in \refeq{d-gauss}
does not contribute, so that \refeq{pi1-error} is bounded by
\begin{align}\lbeq{pi1-errorbd}
\sum_{\substack{t,s,s'\in\Zp\\ 0\leq s\leq s'\leq t}}\frac{O(\beta)}{[1\vee
 (t-s)]^{d/2}}\,\frac{O(\beta)}{(1\vee s')^{d/2}}\leq\sum_{t=0}^\infty
 \frac{O(\beta^2)}{(1\vee t)^{d/2}}\leq O(\beta^2),
\end{align}
where $s,s'$ are the time variables of $\yvec$ and $\zvec$, respectively.
When $\vep<1$, we use the restriction $\omega_i\cap\omega_j=\varnothing$,
$i\ne j$, to extract factors of $q$ with pairwise different arguments,
as in \refeq{pi-premain}, out of the four intersection points
$\ovec,\,\yvec,\,\zvec$ and $\xvec$.  As explained above \refeq{oppi-main},
each pair gives rise to a factor $\vep$, and we obtain a total factor $\vep^4$.
With the help of \refeq{q-gauss}, \refeq{pi1-error} with $\vep<1$ is bounded
by $\vep^{1+4}$ times the left-hand side of \refeq{pi1-errorbd} with the
region of summation being replaced by $\vep\Zp$.  This is further bounded by
$O(\beta^2)\,\vep^2$, since the sum over $t,s,s'\in\vep\Zp$ eats up
a factor $\vep^3$ for the Riemann sum approximation.  This completes the proof of
\refeq{pi1-main}.
\end{proof}

\begin{proof}[Proof of \refeq{intersect1}]
We only consider the case $E^{\rm c}=\Elow^{\uvec}(\omega_1)^{\rm c}$, which is the
event that there is an $\eta\prec\omega_1$ from $\ovec$ to $\xvec$, which must
share at least one step with $\omega_1$, such that $\Elow^{\uvec}(\eta)$ occurs;
the other case $E=\Ehigh^{\uvec}(\omega_2;\omega_1)$ can be estimated in a similar
way.  Let $\omega_3$ be the part of $\eta$ from the point, say $\yvec$,
where $\eta$ starts disagreeing from $\omega_1$ until it hits $\omega_1$
or $\omega_2$ at $\zvec$.  Since $\mP_1$ is a product measure,
\refeq{intersect1} is bounded by
    \[ \vep\sum_{\substack{\xvec,\yvec,\zvec\in\Lambda\\ \zvec\ne\ovec,\xvec}}\,
    \sum_{\substack{\omega_1:\ovec\too\yvec\too\xvec\\ \omega_2:\ovec\too\xvec\\
    \omega_3:\yvec\too\zvec\\ \omega_i\cap\omega_j=\varnothing,~i\ne j}}\big(
    \indic{\zvec\in\omega_1\setminus\{\yvec\}}+\indic{\zvec\in\omega_2}\big)
    \prod_{i=1}^3\mP_1(\omega_i\text{ occupied}). \]

Since the contribution from $\indic{\zvec\in\omega_2}$ is equal to
\refeq{pi1-error}, we only need to investigate the contribution due to the
other indicator $\indic{\zvec\in\omega_1\setminus\{\yvec\}}$.  We again
discuss the case $\vep=1$ first, and then adapt the argument to the case
$\vep<1$, as done below \refeq{pi1-errorbd}.  When $\vep=1$, we ignore
the restriction $\omega_i\cap\omega_j=\varnothing$, $i\ne j$, and apply
\refeq{d-gauss} to the probability of $\omega_2$ and $\omega_3$ being
occupied.  By denoting the time variables of $\yvec$ and
$\zvec$ by $s$ and $s'$ respectively, the contribution from
$\indic{\zvec\in\omega_1\setminus\{\yvec\}}$ is bounded by
\begin{align}
\sum_{\substack{t,s,s'\in\Zp\\ 0\leq s<s'\leq t}}\frac{O(\beta)}{(1\vee t)
 ^{d/2}}\,\frac{O(\beta)}{[1\vee(s'-s)]^{d/2}}\leq\sum_{t=0}^\infty\frac{O
 (\beta^2)}{(1\vee t)^{(d-2)/2}}\leq O(\beta^2).
\end{align}
When $\vep<1$, we use the restriction $\omega_i\cap\omega_j=\varnothing$,
$i\ne j$, along each of the four intersection points and obtain the eight
factors of $q$ with pairwise different arguments.  Following the argument
below \refeq{pi1-errorbd}, we obtain the desired bound $O(\beta^2)\,\vep^2$.
This completes the proof of \refeq{intersect1}.
\end{proof}

\begin{proof}[Proof of \refeq{intersect2}]
Since $\omega_1\cap\omega_2\ne\varnothing$, there is a sequence of bonds
$b_1,\dots,b_n$ such that $\omega_1$ and $\omega_2$ meet for the first time
at $\bb_1$, share $b_1,\dots,b_n$, and split at $\tb_n$ ($\omega_1$ and
$\omega_2$ may share a bond again after $\tb_n$).  This means that, together
with $q(u)\,q(v)\,q(x-y)\,q(x-z)$ in \refeq{intersect2}, the left-hand side
of \refeq{intersect2} is bounded by the convolution of two non-vanishing
bubbles and $\prod_{i=1}^nq(w_i)^2$, where each $w_i$ is the spatial
component of $\tb_i-\bb_i$.  Using \refeq{q-gauss} and
$\sum_wq(w)^2\leq\|q\|_\infty$, we can bound
\refeq{intersect2} by
\begin{align}
\vep\sum_{\substack{t,s,s'\in\vep\Zp\\ \vep<s<s'<t-\vep}}\frac{O(\beta)\,
 \vep^2}{(1\vee s)^{d/2}}~\big[(1-\vep)\vee(\vep\|D\|_\infty)\big]^{\frac
 {s'-s}\vep}\,\frac{O(\beta)\,\vep^2}{[1\vee(t-s')]^{d/2}}\leq O(\beta^2)
 \,\vep^2,
\end{align}
where, as before, $\vep^3$ is used up for the Riemann sum approximation.
The above estimate can be improved to $O(\beta^3)$ for oriented percolation,
using \refeq{Dprop}.  This completes the proof of \refeq{intersect2}.
\end{proof}




\subsection{Asymptotics of $\pcperc$}\label{ss:perc}
In this section, we compute the asymptotics of the critical point for
(unoriented) percolation.  We follow the strategy in Section \ref{ss:op} as
closely as possible.  However, there are a number of changes due to the fact
that we have less control of the lace expansion coefficients.  For example,
the bounds on the derivative of $\hat\Pi_p$ with respect to $p$ are not
available in the literature, even though in the unpublished manuscript
\cite{h03}, this derivative is computed.  To make this paper self-contained,
we avoid the use of the derivative, which causes changes in the proof.

We start with some notation. Let
\begin{align}\lbeq{Tdef}
T_p=\sup_{x\in \Z^d} (pD*\tau_p^{*3})(x),&&
T_p'=\sup_{x\in \Z^d} \tau_p^{*3}(x).
\end{align}
We will use the following bounds:

\begin{lemma} Fix $d>6$. For $L$ sufficiently large, and all $p\leq \pc$,
\label{lem-T}
\begin{align}
T_p\leq C\beta,&& T_p'\leq1+C\beta.\lbeq{Tbd}
\end{align}
\end{lemma}
We will defer the proof of Lemma \ref{lem-T} to the end of this section.

To compute the asymptotics of $\hat\Pi_p$, we use \refeq{perc-lace}
and the bound (cf., \cite[Proposition 4.1]{BCHSS04b})
\begin{align}\lbeq{piNbd}
\hat{\pi}^\smallsup{N}_p \leq T_p'(2T_pT_p')^{N\vee 1}.
\end{align}
Note that Lemma \ref{lem-T} together with \refeq{perc-pc} and \refeq{piNbd}
immediately imply
\begin{align}
\pc=1+O(\beta).\lbeq{pcasyone}
\end{align}

We now start the proof to improve \refeq{pcasyone} one term further.
Together with Lemma \ref{lem-T}, \refeq{piNbd} proves that the contribution
to $\sum_{N=2}^{\infty}\hat\pi_{\pc}^{\sss(N)}$ is $O(\beta^2)$.  Thus, we
are left to compute $\hat\pi_{\pc}^{\sss(0)}$ and $\hat\pi_{\pc}^{\sss(1)}$.
The goal of this section is to prove
\begin{align}\lbeq{pi01asy}
\hat{\pi}_{\pc}^{\scscst(0)}=\frac12\sum_{n=3}^{\infty}(n-1)D^{*n}(o)
 +O(\beta^2),&&
\hat{\pi}_{\pc}^{\scscst(1)}=D^{*2}(o)+\sum_{n=3}^{\infty}nD^{*n}(o)
 +O(\beta^2).
\end{align}
Using \refeq{perc-pc} and \refeq{pi01asy}, we arrive at \refeq{perc-critpt}.
Thus, we are left to prove \refeq{pi01asy}.

We again investigate $\hat\pi_{\pc}^{\sss(0)}$ and $\hat\pi_{\pc}^{\sss(1)}$
separately.  First, we compute $\hat\pi_{\pc}^{\sss(0)}$. For percolation,
we denote by $\{w\udb x\}$ the event that $w$ is doubly connected to $x$.
By definition \cite{HS90a}, 
$\hat\pi_p^{\sss(0)}=\sum_{x\in\Zd}\pi_p^{\sss(0)}(x)$, where
\begin{align}
\pi_p^{\sss(0)}(x)=\mP_p(o\udb x)-\delta_{o,x}.
\end{align}
We wish to use {\it Russo's formula} (see, e.g.,
\cite{g99}) to prove that $\hat \pi_{\pc}^{\scscst(0)}=
\hat \pi_1^{\scscst(0)}+O(\beta^2)$. However, Russo's formula
is restricted to events that only depend on a {\it finite} number of bonds,
so that we will first show that Russo's formula may be applied to
$\pi_p^{\scscst(0)}(x)$.

Let $B_\ell=\{x\in\Zd:|x|\leq\ell\}$.  We note that, since
$\hat\pi_p^{\sss(0)}$ is finite for any $p\leq\pc$, there is an $r<\infty$
such that $\sum_{x\notin B_r}\pi_p^{\sss(0)}(x)=O(\beta^2)$ for any
$p\leq\pc$.  In fact, using the {\it BK inequality} (see, e.g., \cite{g99})
and the bound $\tau_{\pc}(x)\leq K|x|^{2-d}$ for $x\ne o$
\cite[Proposition~2.2]{HHS01a}\footnote{In \cite[Proposition~2.2]{HHS01a},
$K$ is of order $O(L^{-2+\epsilon})$ with an arbitrarily small number
$\epsilon>0$, and thus is small when $L$ is large.  Here, we do not care
about the dependence of $K$ on $L$, and will take $K=O(1)$.}\!, we have
    \begin{align}\lbeq{finitesum}
    \sum_{x\notin B_r}\pi_p^{\sss(0)}(x)\leq\sum_{x\notin B_r}\tau_p(x)^2
    \leq\sum_{x\notin B_r}\tau_{\pc}(x)^2\leq c\sum_{\ell>r}\ell^{(d-1)
    +2(2-d)}=O(r^{4-d})=O(\beta^2),
    \end{align}
where we assume $r=O(L^{2d/(d-4)})$.  Let $\{E$ in $B_{\sR}\}$ be the set of
bond configurations whose restriction on bonds $\{u,v\}$ with $u,v\in B_{\sR}$
are in $E$.  Similarly to \refeq{finitesum}, if $R=O(L^{2d/(d-6)})$, then
for any $p\leq\pc$ we have\footnote{The event
$\{o\udb x\}\setminus\{o\udb x\text{ in }B_{\sR}\}$ implies the existence of
$y\notin B_{\sR}$ such that $o\cntd x$, $x\cntd y$ and $y\cntd o$ occur
disjointly.  Therefore, by the BK inequality and Propositions~1.7(i)
and 2.2 of \cite{HHS01a}, we obtain
    \[ \sum_{x\in\Zd}\mP_p\big(\{o\udb x\}\setminus\{o\udb x\text{ in }B_{\sR}\}
    \big)\leq\sum_{\substack{x\in\Zd\\ y\notin B_{\sR}}}\tau_p(x)~\tau_p(y-x)~
    \tau_p(y)\leq c\sum_{y\notin B_{\sR}}|y|^{(4-d)+(2-d)}=O(R^{6-d}). \]
    }
    \begin{align}\lbeq{finitecont}
    \sum_x\mP_p\big(\{o\udb x\}\setminus\{o\udb x\text{ in }B_{\sR}\}\big)
    \leq O(\beta^2).
    \end{align}
By \refeq{finitesum}--\refeq{finitecont} and the mean-value theorem,
there is a $p\in(1,\pc)$ such that
    \begin{align}
    \hat\pi_{\pc}^{\sss(0)}&=\sum_{x\in B_r}\pi_{\pc}^{\sss(0)}(x)+O(\beta^2)=
    \sum_{x\in B_r\setminus\{o\}}\mP_{\pc}(o\udb x\text{ in }B_{\sR})+O(\beta^2)\nn\\
    &=\sum_{x\in B_r\setminus\{o\}}\mP_1(o\udb x\text{ in }B_{\sR})+(\pc-1)\sum_{x\in
    B_r\setminus\{o\}}\partial_p\mP_p(o\udb x\text{ in }B_{\sR})+O(\beta^2)\nn\\
    &=\hat\pi_1^{\sss(0)}+(\pc-1)\sum_{x\in B_r\setminus\{o\}}\partial_p\mP_p(o
    \udb x\text{ in }B_{\sR})+O(\beta^2).\lbeq{mvt}
    \end{align}
We will later identify $\hat\pi_1^{\sss(0)}$, and first show that the second
term is $O(\beta^2)$.  Since the event $\{o\udb x$ in $B_{\sR}\}$ depends only
on finitely many bonds, we are now allowed to apply Russo's formula to obtain
    \begin{align}
    \sum_{x\in B_r\setminus\{o\}}\partial_p\mP_p(o\udb x\text{ in }B_{\sR})=
    \sum_{x\in B_r\setminus\{o\}}\,\sum_{(u,v)}D(v-u)~\mP_p\big((u,v)
    \text{ pivotal for }\{o\udb x\text{ in }B_{\sR}\}\big),\lbeq{russo}
    \end{align}
where the factor $D(v-u)$ arises from the derivative of the bond occupation
probability of $\{u,v\}$ with respect to $p$, and where a bond is pivotal for
$o\udb x$ when $o\udb x$ in the (possibly modified) configuration where
the bond is made occupied, while $o\udb x$ does not occur in the
(possibly modified) configuration where the bond is made occupied.

Since $\pc=1+O(\beta)$, and since, by the BK inequality, \refeq{russo} is bounded by
    \begin{align}
    &\sum_{x,(u,v)}D(v-u)~\mP_p\big(\{o\cntd u\}\circ\{v
    \cntd x\}\circ\{o\cntd x\}\big)\nn\\
    &\qquad\leq\sum_{x,(u,v)}D(v-u)~\tau_p(u)\,\tau_p(x-v)\,\tau_p
    (x)\leq p^{-1}T_p\leq T_p,
    \end{align}
so that the second term in \refeq{mvt} is $O(\beta^2)$.  We are left to
analyse the first term $\hat\pi_1^{\sss(0)}$.  We follow the strategy
around \refeq{pi0-eventreexpr}, but the details change somewhat.


Let $\Scal_x$ denote all self-avoiding paths from $o$ to $x$, and order
the elements in $\Scal_x$ in an arbitrary way.  Then we can write
\begin{align}\lbeq{pi-reexprup}
\hat\pi_1^{\sss(0)}=\sum_{x\ne o}~\sum_{\substack{\omega_1,\omega_2\in
 \Scal_x\\ \omega_1\cap\omega_2=\varnothing\\ \omega_1\prec\omega_2}}\mP_1
 \big(\omega_1,\omega_2\text{ occupied};\,\Elow(\omega_1)\cap\Ehigh(\omega_2;
 \omega_1)\big),
\end{align}
where $\Elow(\omega_1)$ and $\Ehigh(\omega_2;\omega_1)$ were defined between
\refeq{pi0-main} and \refeq{pi0-eventreexpr}.  In words, the event
$\Elow(\omega_1)$ holds when $\omega_1$ is the lowest occupied self-avoiding
walk path from $o$ to $x$ such that there is an occupied bond disjoint path
from $o$ to $x$.  The event $\Ehigh(\omega_2;\omega_1)$ holds when $\omega_2$
is the highest occupied self-avoiding walk path from $o$ to $x$ that is bond
disjoint from $\omega_1$.  Since $\mP_1$ is a product measure, if we ignore
$\Elow(\omega_1)\cap\Ehigh(\omega_2;\omega_1)$, we obtain
\begin{align}\lbeq{pi-premainup}
\sum_{\substack{\omega_1,\omega_2\in\Scal_x\\ \omega_1\cap\omega_2
 =\varnothing\\ \omega_1\prec\omega_2}}\mP_1(\omega_1,\omega_2
 \text{ occupied})=\sum_{\substack{\omega_1,\omega_2\in\Scal_x\\
 \omega_1\cap\omega_2=\varnothing\\ \omega_1\prec\omega_2}}\mP_1
 (\omega_1\text{ occupied})~\mP_1(\omega_2\text{ occupied}).
\end{align}
We can then follow the rest of the argument between \refeq{pi-premain}
and \refeq{oppi-main} to arrive at
\begin{align}
\hat\pi_1^{\sss(0)}=\frac12\sum_{x\ne o}\sum_{\substack {\omega_1,\omega_2
 \in\Scal_x\\ |\omega_1|+|\omega_2|\geq3}}W_1(\omega_1)\;W_1(\omega_2)+
 O(\beta^2),\lbeq{percpi0-main}
\end{align}
where we recall the definition of $W_p(\omega)$ in \refeq{weight-def}.  Here
the factor $1/2$ has the same origin as the one in \refeq{decomp}, and the
restriction that $|\omega_1|+|\omega_2|\geq3$ is due to the fact that the
smallest cycle in percolation has length 3.  In \refeq{percpi0-main}, each
$\omega_j$ is a self-avoiding path from $o$ to $x$.  However, as estimated
in Section~\ref{s:saw}, the contribution in which $\omega_1$ or $\omega_2$
has a self-intersection is $O(\beta^2)$.  Therefore, we can remove the
self-avoidance constraint in \refeq{percpi0-main}.  Performing the sum
over $x\ne o$ and writing $\omega=(\omega_1,\omega_2)$, which is a random
walk path starting and ending at $o$ of length at least 3, we obtain
\begin{align}
\hat\pi_1^{\sss(0)}=\frac12\sum_{\substack{\omega:o\too o\\ |\omega|\geq3}}
 (|\omega|-1)\,W_1(\omega)+O(\beta^2)=\frac12\sum_{n=3}^\infty(n-1)\,D^{*n}
 (o)+O(\beta^2),
\end{align}
where $|\omega|-1=n-1$ is the number of vertices along $\omega$, excluding
the starting and ending point of $\omega$.
This completes the computation of the leading asymptotics of
$\hat\pi_{\pc}^{\sss(0)}$.

We next derive the asymptotics of $\hat{\pi}_{\pc}^{\scscst(1)}$, following
the strategy in \cite{HS04a, HS04b}, where the first three coefficients of
the asymptotic expansion into powers of $(2d)^{-1}$ of the critical value
$\pc$ for nearest-neighbour percolation were computed.  The details of the
argument are changed considerably compared to \cite{HS04a, HS04b}.  Indeed,
since we are only interested in the leading order term, while in \cite{HS04a}
the first three coefficients are computed, many terms that need explicit
computation in \cite{HS04a, HS04b} will be error terms for us.  On the other
hand, since in \cite{HS04a, HS04b} the asymptotics in nearest-neighbour
models for large dimensions are considered, long loops lead to error term in
\cite{HS04a, HS04b}, whereas they contribute to the leading asymptotics here.
We follow the proof in \cite[Section 4.2]{HS04a} as closely and as long as
possible, and indicate where the argument diverges.

To define $\hat{\pi}_p^{\smallsup{1}}$, we need the following definitions.
Given a bond configuration and $A\subseteq\Zd$, we recall that $x$ and $y$ are
\emph{connected through} $A$, and write $x\ctx{\;A\;}y$, if every occupied
path connecting $x$ to $y$ has at least one bond with an endpoint in $A$.
As defined below \refeq{op-pi0}, the directed bond $(u,v)$ is said to be
{\em pivotal} for $x\cntd y$, if $x\cntd u$ and $v\cntd y$ occur, and if
$x\cntd y$ occurs only when $\{u,v\}$ is set occupied.  (Note that there is
a distinction between the events $\{(u,v)$ is pivotal for $x\cntd y\}$ and
$\{(v,u)$ is pivotal for $x\cntd y\}=\{(u,v)$ is pivotal for $y\cntd x\}$.)
Let
\begin{align}\lbeq{E'def}
E'(v,x;A)=\{v\ctx{\;A\;}x\}\cap\big\{\nexists(u',v')
 \text{ occupied \& pivotal for }v\cntd x\text{ s.t. }v\ctx{\;A\;}u'\big\}.
\end{align}
Then, by definition \cite{HS90a},
\begin{align}\lbeq{Pi1defa}
\hat\pi_p^{\sss(1)}=\sum_x\sum_{(u,v)}p\,D(v-u)~\mE_{\sss0}\Big[\ind{o\udb
 u}~\mP_{\sss1}\big(E'(v,x;\tilde C^{\sss(u,v)}_{\sss0}(o))\big)\Big],
\end{align}
where the sum over $(u,v)$ is a sum over directed bonds.  On the right-hand
side, we use subscripts to identify the different expectations.  Thus, the
subscripts do {\it not} refer to the percolation parameter $p$.  The cluster
$\tilde C^{\sss(u,v)}_{\sss0}(o)$ appearing on the right hand side of
\refeq{Pi1defa} is random with respect to the expectation $\mE_{\sss0}$,
but $\tilde C^{\sss(u,v)}_{\sss0}(o)$ should be regarded as a \emph{fixed} set
inside the probability $\mP_{\sss1}$.  The latter introduces a second
percolation model which depends on the original percolation model via the
set $\tilde C^{\sss(u,v)}_{\sss0}(o)$.  We refer to the bond configuration
corresponding to the $j^{\rm th}$-expectation as the ``level-$j$''
configuration.

By \refeq{piNbd},
\begin{align}\lbeq{Pi1Tbd}
0\leq\hat\pi_p^\smallsup{1}\leq 2T_p'T_p T_p'.
\end{align}
We will use refinements of this bound in the following.

We first claim that the contribution to \refeq{Pi1defa} due to $u\ne o$
is an error term of order $O(\beta^2)$.  Indeed, if $u\neq o$ then at
level-0 the origin is in a cycle of length at least 3.  Standard
diagrammatic estimates then allow for the replacement in \refeq{Pi1Tbd}
of a factor $T_p'$ by a constant multiple of $T_p$.  This improves the
bound \refeq{Pi1Tbd} from $O(\beta)$ to $O(\beta^2)$, by \refeq{Tbd}.

We are left with the contribution to \refeq{Pi1defa} due to $u=o$, namely
\begin{align}\lbeq{Pi1defu0}
\sum_{x,v}p\,D(v)~\mE_{\sss0}\big[\mP_{\sss1}\big(E'(v,x;\tilde C^{\sss
 (o,v)}_{\sss0}(o))\big)\big].
\end{align}
If $ x\notin\tilde C^{\sss(o,v)}_{\sss0}(o)$, then to obtain a non-zero
contribution to $\mP_{\sss1}(E'(v,x;\tilde C^{\sss(o,v)}_{\sss0}(o)))$,
$x$ must be in an occupied cycle of length at least 3, in level-1
(in the language of \cite[Section~3]{BCHSS04b}, the sausage containing
$x$ must consist of a cycle containing both $x$ and an endpoint of the
last pivotal bond for the connection from $o$ to $x$).  In this case,
in \refeq{Pi1Tbd}, we may again replace a factor $T_p'$ by a constant
multiple of $T_p$, and again this contribution is $O(\beta^2)$.  We are
left to consider
\begin{align}\lbeq{Pi1defux}
\sum_{x,v}p\,D(v)~\mE_{\sss0}\Big[\ind{x\in\tilde C^{\sss(o,v)}_0(o)}
 ~\mP_{\sss1}\big(E'(v,x;\tilde C^{\sss(o,v)}_{\sss0}(o))\big)\Big].
\end{align}
This is as far as the analogy with the argument in \cite[Section~4.2]{HS04a}
goes.  We now need to adapt the proof there to compute the asymptotics of
$\hat\pi^{\sss(1)}_{\pc}$ when $L\to\infty$.

If $x\in\tilde C^{\sss(o,v)}_{\sss0}(o)$, and if $v\cntd x$, then
$v\ctx{\tilde C^{\sss(o,v)}_{\sss0}(o)}x$.  We next claim that the
intersection with the second event in \refeq{E'def} leads to an error
term.  We write
\begin{align*}
&I_{\sss0}[x\in\tilde C^{\sss(o,v)}_{\sss0}(o)]~I_{\sss1}[E'(v,x;\tilde
 C^{\sss(o,v)}_{\sss0}(o))]\\[5pt]
&\qquad=I_{\sss0}[x\in\tilde C^{\sss(o,v)}_{\sss0}(o)]~I_{\sss1}[v\cntd
 x]\nn\\
&\qquad\qquad\times\Big(1-I_{\sss1}\big[\exists(u',v')
 \text{ occupied \& pivotal for }v\cntd x\text{ s.t. }
 v\ctx{\tilde C^{\sss(o,v)}_{\sss0}(o)}u'\big]\Big),
\end{align*}

\noindent
where we write $I_{\sss0}$ and $I_{\sss1}$ for the indicator functions on
levels 0 and 1, respectively.  The latter term can be bounded by
\begin{align}
\sum_{(u',v')}\sum_zI_{\sss0}[z\in\tilde C^{\sss(o,v)}_{\sss0}(o)]~I_{\sss1}
 \big[\{v\cntd z\}\circ\{z\cntd u'\}\circ\{(u',v')\text{ occupied}\}\circ
 \{v'\cntd x\}\big],
\end{align}
which, using the BK inequality, yields a bound of the form
\begin{align}
\sum_{x,v,z}\sum_{(u',v')}p\,D(v)~\mP_{\sss0}(o\cntd x,~o\cntd z)~\tau_p
 (z-v)~\tau_p(u'-z)~p\,D(v'-u')~\tau_p(x-v').
\end{align}
By the tree-graph inequality \cite{an84}
\begin{align}
\mP_{\sss0}(o\cntd x,~o\cntd z)\leq\sum_y\tau_p(y)~\tau_p(x-y)~\tau_p(z-y),
\end{align}
so that we end up with
\begin{align}
\sum_{x,z,y}(p\,D*\tau_p)(y)~\tau_p(x-y)~\tau_p(z-y)~\tau_p(z-v)~(\tau_p*p
 \,D*\tau_p)(x-z)\leq T_p^2=O(\beta^2),
\end{align}
which indeed is an error term.  Thus, using the identity
\begin{align}
\{x\in\tilde C^{\sss(o,v)}_{\sss0}(o)\}=\{o\cntd x
 \text{ without using }(o,v)\},
\end{align}
we end up with
\begin{align}\lbeq{pi1asya}
\hat\pi_p^{\sss(1)}=\sum_{v,x}p\,D(v)~\tau_p^{\sss(o,v)}(x)~\tau_p(x-v)
 +O(\beta^2),
\end{align}
where
\begin{align}\lbeq{tau(o,v)def}
\tau_p^{\sss(o,v)}(x)=\mP(o\cntd x\text{ without using }(o,v)).
\end{align}
Note that we can think of $\tau_p^{\sss(o,v)}(x)$ as the two-point function
on $\Zd$, where the bond $(o,v)$ is removed. We will denote the resulting
graph with vertex set $\Z^d$ and edge set
$\big\{\{x,y\}:x,y\in\Zd,\{x,y\}\ne\{o,v\}\big\}$ by $\zd_{\sss(o,v)}$,
so that $\tau_p^{\sss(o,v)}(x)$ is the two-point function on
$\zd_{\sss(o,v)}$.  We will use this observation to compute
$\tau_p^{\sss(o,v)}(x)$.

We investigate the main term in the right-hand side of \refeq{pi1asya}
further.
%
%
%
Russo's formula, together with the BK inequality, yields that
\begin{align}
\partial_p\tau_p(x)&=\sum_{(y,z)}D(z-y)~\mP((y,z)\text{ pivotal for }
 o\cntd x)\leq(\tau_p*D*\tau_p)(x),\\
\partial_p\tau_p^{\smallsup{o,v}}(x)&=\sum_{(y,z)}D(z-y)~\mP((y,z)
 \text{ pivotal for }o\cntd x\text{ in }\zd_{\smallsup{o,v}})
 \leq(\tau_p*D*\tau_p)(x).
\end{align}
Therefore, we obtain that for $p=\pc$,
\begin{align}
\hat{\pi}_{\pc}^\smallsup{1} &=\sum_{x,v}\pc D(v) \tau_{\pc}^{\smallsup{o,v}}(x)\tau_{\pc}(x-v)
+O(\beta^2)\nn\\
&=\sum_{x,v}D(v) \tau_{1}^{\smallsup{o,v}}(x)\tau_{1}(x-v)
+O((\pc-1)T_{\pc})+O(\beta^2)\nn\\
&=\sum_{x,v}D(v) \tau_{1}^{\smallsup{o,v}}(x)\tau_{1}(x-v)
+O(\beta^2),
\end{align}
since $\pc=1+O(\beta)$. Furthermore, an argument similar to the one for
$\hat{\pi}_p^{\smallsup{0}}$ shows that
\begin{align}
\tau_1(x)&=G(x)+O\big((G*g*G)(x)\big),\lbeq{taup1}\\
\tau_1^{\sss(o,v)}(x)&=G(x)\,(1-\delta_{v,x})+\delta_{v,x}(D^{*2}*G)(x)
 +O\big((G*g*G)(x)\big)+O\big(D*(G-\delta_o)(x)\big).
\lbeq{taupv1}
\end{align}
where we recall $G(x)=\sum_{n=0}^\infty D^{*n}(x)$ and define
\begin{align}
g(x)=G(x)(D*G)(x).
\end{align}
We will prove \refeq{taup1}--\refeq{taupv1} in full detail below, and first
complete the proof subject to \refeq{taup1}--\refeq{taupv1}.  Using
\refeq{taup1}--\refeq{taupv1}, together with the fact that for $u\neq o$,
we have $G(u)=(D*G)(u)$, we end up with
\begin{align}\lbeq{pi1asyconc}
\hat\pi_{\pc}^\smallsup{1}&=\sum_{x\neq v}G(x)\,D(v)\,G(x-v)+\sum_x(D^{*2}*G)
 (x)\,D(x)+O(\beta^2)+O((D*G^{*3}*g)(o))\nn\\
&=\sum_{n=2}^\infty(n-1)D^{*n}(o)+\sum_{n=3}^\infty D^{*n}(o)+O(\beta^2)
 =D^{*2}(o) +\sum_{n=3}^{\infty} n D^{*n}(o)+O(\beta^2),
\end{align}
where we use
\begin{align}
(D*G^{*3}*g)(o)\leq\|D*G^{*3}\|_\infty\,\|g\|_1\leq O(\beta^2),
\end{align}
for $d>6$, by \refeq{d-gauss}.

This completes the proof subject to \refeq{taup1}--\refeq{taupv1} and
Lemma~\ref{lem-T}.
\qed

\newcommand{\Flow}{F_{\succ}}

\begin{proof}[Proof of \refeq{taup1}--\refeq{taupv1}]
We start by proving \refeq{taup1}, and then adapt the argument to prove
\refeq{taupv1}.  To see \refeq{taup1}, we recall the arbitrary ordering of
the elements in $\Scal_x$ introduced above \refeq{pi-reexprup}.  Then we
have that
\eq
\tau_1(x)=\sum_{\omega\in\Scal_x}\mP_p(\omega\text{ occupied};\Flow(\omega)),
\en
where $\Flow(\omega)$ is the event that $\omega$ is the lowest occupied path
in $\Scal_x$.  Thus, we can write
\eq
\tau_1(x) = \sum_{\omega\in \Scal_x} \mathbb{P}_p(\omega \text{ occupied}) -
\sum_{\omega\in \Scal_x} \mathbb{P}_p(\omega \text{ occupied}; \Flow(\omega)^c).
\lbeq{tau1rew}
\en
The former term equals
\eq
\lbeq{sawpaths}
\delta_{o,x}+(1-\delta_{o,x})\sum_{\omega\in \Scal_x} \prod_{i=0}^{|\omega|-1}D(\omega(i+1)-\omega(i)).
\en
Clearly, by using inclusion-exclusion on the fact that $\omega$ is self-avoiding,
as in \refeq{saw-pi1exp}, \refeq{sawpaths} equals
\eq
G(x)+O((G*G)(x)(G(o)-1)),
\en
which is a contribution to the error in \refeq{taup1} when we note that
$G(o)-1=(D*G)(o)$. Similarly, the second term in \refeq{tau1rew} is bounded
by $O\big((G*g*G)(x)\big)$ using the fact that there must exist a $u\in\Zd$
such that there exist bond disjoint occupied paths from $o$ to $u$, two
occupied paths from $u$ to $v$ (of which at least one is non-vanishing) and
one from $v$ to $x$. Thus, by the BK inequality, this term is bounded by
\[ \sum_{u,v}G(u)\,G(v-u)\,(D*G)(v-u)\,G(x-v)=(G*g*G)(x). \]

The proof of \refeq{taupv1} follows the same ideas. In
\refeq{tau1rew} and \refeq{sawpaths}, we only need to sum over
self-avoiding walk paths that do not use the bond
$(o,v)$. When $x=v$, this means that $|\omega|\geq 2$, so that we obtain
\eq
\tau_1^{\smallsup{o,v}}(v)=(D^{*2}*G)(v)+O((G*g*G)(v)) + O(D(v)(G(o)-1)).
\en
When $x\neq v$, we can use inclusion-exclusion on the fact that the bond
$(o,v)$ is not used, and obtain
\eq
\tau_1^{\smallsup{o,v}}(x) = \tau_1(x)+O(D(v) G(x-v)),
\en
and then use \refeq{taup1}.
\end{proof}

\begin{proof}[Proof of Lemma \ref{lem-T}]
We use \cite[(5.20)]{HS90a}, which states that uniformly in $p\in[1,\pc)$
and for $L$ large enough
\begin{align}\lbeq{IRbd}
\hat\tau_p(k)\leq\frac{1+o(1)}{1-\hat{D}(k)},
\end{align}
where $o(1)$ tends to 0 when $L\rightarrow \infty$. 
%
We also use the standard bound (see e.g.\ \cite{BCHSS04b}) that for $x\neq 0$,
\eq
\tau_p(x) \leq (pD*\tau_p)(x).
\lbeq{tauDtau}
\en
We then follow the proof as in \cite{BCHSS04b}.  For $T_p$, we fix $x$ and
extract the term in \refeq{Tdef} due to the case where every argument of
$\tau_p$ is $o$, which is $pD(x)\leq p\,C\beta$ (see \refeq{Dprop}).
This gives
\begin{align}
T_p(x)\leq p\,C\beta+p\sum_{(u,y,z)\ne(x,o,o)}\tau_p(y)\,\tau_p(z-y)\,D(u)\,
 \tau_p(x+z-u).
\end{align}
Therefore, by \refeq{tauDtau},
\begin{align}
T_p\leq p\,C\beta+3p^2\sup_x(D^{*2}*\tau_p^{*3})(x),
\end{align}
where the factor 3 comes from the 3 factors $\tau_p$ whose argument can
differ from $o$.  In terms of the Fourier transform, this gives
\begin{align}
T_p\leq p\,C\beta +3p^2\,\sup_x\int_{\Box_\pi}\frac{d^d k}{(2\pi)^d}~\hat
 D(k)^2\,\hat\tau_p(k)^3\,e^{-ik\cdot x}=p\,C\beta+3p^2\int_{\Box_\pi}
 \frac{d^dk}{(2\pi)^d}~\hat D(k)^2\,\hat\tau_p(k)^3,
\end{align}
where $\Box_\pi=[-\pi,\pi]^d$ and we use $\hat{\tau}_p(k)\geq 0$
\cite{an84}.  For $L\gg1$, by \refeq{IRbd} we obtain
    \begin{align}\lbeq{Tpintbd}
    T_p\leq p\,C\beta+4p^2\int_{\Box_\pi}\frac{d^d k}{(2\pi)^d}~\frac{\hat
    D(k)^2}{[1-\hat D(k)]^3}.
    \end{align}
Let $\hat B_{1/L}=\{k\in\Box_\pi:cL^2|k|^2\leq\eta\}$.  Using \refeq{Dprop},
we estimate the contribution to the integral in \refeq{Tpintbd} from
$k\in\Box_\pi\setminus\hat B_{1/L}$ by
    \begin{align}
    \int_{\Box_\pi\setminus\hat B_{1/L}}\frac{d^d k}{(2\pi)^d}~\frac{\hat D
    (k)^2}{[1-\hat D(k)]^3}\leq\eta^{-3}\int_{\Box_\pi}\frac{d^d k}{(2\pi)^d}
    ~\hat D(k)^2=O(\beta).
    \end{align}
On the other hand, the contribution from $k\in\hat B_{1/L}$ is, again using
\refeq{Dprop}, bounded by
    \begin{align}
    \int_{\hat B_{1/L}}\frac{d^d k}{(2\pi)^d}~\frac{\hat{D}(k)^2}{[1-\hat{D}(k)]
    ^3}\leq\int_{\hat B_{1/L}}\frac{d^d k}{(2\pi)^d}~(cL^2|k|^2)^{-3}=O(\beta).
    \end{align}
This proves the bound on $T_p$.

The bound on $T_p'$ is a consequence of $T_p' \leq 1+3T_p$.  Here the term
1 is due to the contribution to \refeq{Tdef} where the arguments of the
three factors of $\tau_p$ in $T_p'$ in \refeq{Tdef} are equal to $o$.
If at least one of these arguments is nonzero, then we can use
\refeq{tauDtau} for the corresponding two-point function.
\end{proof}

\appendix
\section{Bounds on $D^{*n}(x)$ and $q^{*n}(x)$}\label{s:lclt}

In this appendix, we prove \refeq{q-gauss} for any $\vep\in (0,1]$
assuming \refeq{Dprop}. The inequality \refeq{d-gauss} follows
by taking $\vep=1$.

First, we note that
\begin{align}
q^{*n}(x)=(1-\vep)^n\delta_{o,x}+\sum_{j=0}^{n-1}(1-\vep)^{n-1-j}
 \vep(D*q^{*j})(x),\lbeq{q-dec}
\end{align}
where we suppose that the empty sum equals zero.  When
$n\leq N\equiv\vep^{-1}$, we use \refeq{q-dec} to obtain
    \begin{align}\lbeq{qbd-small}
    q^{*n}(x)\leq(1-\vep)^n\delta_{o,x}+\sum_{l=0}^{n-1}(1-\vep)^l
    \vep\,\|D\|_\infty\leq(1-\vep)^n\delta_{o,x}+O(\beta)\leq (1-\vep)^n\delta_{o,x}
    +\frac{O(\beta)}{[1\vee (n\vep)]^{d/2}},
    \end{align}
as required. On the other hand, when $n>N$, we use
    \begin{align}
    q^{*n}(x)=(1-\vep)^n\delta_{o,x}+(1-\vep)^{n-1}n\vep D(x)+\sum_{j=0}^{n-2}
    (n-1-j)(1-\vep)^{n-2-j}\vep^2(D^{*2}*q^{*j})(x),\lbeq{q-dec2}
    \end{align}
which is obtained by substituting \refeq{q-dec} into $q^{*j}$ in the
right-hand side of \refeq{q-dec}.  Since the
second term in the right-hand side is bounded by $O(\beta)n\vep e^{-n\vep}
\leq O(\beta)[1\vee (n\vep)]^{-d/2}$, it suffices to investigate the third term.
Let $S_1$ be the sum over $j<N$, and let $S_2$ be the remaining sum, i.e.,
    \begin{align}
    S_1&=\sum_{0\leq j<N}(n-1-j)(1-\vep)^{n-2-j}\vep^2(D^{*2}*q^{*j})(x),\\
    S_2&=\sum_{N\leq j\leq n-2}(n-1-j)(1-\vep)^{n-2-j}\vep^2(D^{*2}*q^{*j})(x).
    \lbeq{S2}
    \end{align}
For $S_1$, we use
\refeq{qbd-small} to obtain
    \begin{align}\lbeq{S1-bd}
    S_1
    \leq(1-\vep)^{n-2}\vep^2D^{*2}(o)\sum_{l=n-N}^{n-1}l+O(\beta)\,\vep^2
    \sum_{l=n-N}^{n-1}l\,(1-\vep)^{l-1}\leq O(\beta)\,(n\vep)^2e^{-n\vep}.
    \end{align}
For $S_2$, we recall the definition $\hat B_{1/L}=\{k\in\Box_\pi:cL^2|k|^2\leq\eta\}$
below \refeq{Tpintbd}, and let $\hat B_{1/L}^+=\{k\in\hat B_{1/L}:\hat q(k)\geq0\}$.
Note that
    \begin{align}
    (D^{*2}*q^{*j})(x)\leq\int_{\hat B_{1/L}^+}\frac{d^dk}{(2\pi)^d}~\hat
    q(k)^j+\int_{\Box_\pi\setminus\hat B_{1/L}^+}\frac{d^dk}{(2\pi)^d}~
    \hat D(k)^2|\hat q(k)|^j.\lbeq{D*q-prebd}
    \end{align}
Recall \refeq{Dprop}.  For the first integral, we use
    \begin{align}
    \hat q(k)=1-\vep[1-\hat D(k)]\leq e^{-\vep[1-\hat D(k)]}\leq e^{-c\vep
    L^2|k|^2},\lbeq{q-bd}
    \end{align}
while for the second integral in \refeq{D*q-prebd}, we use, noting that without loss of generality,
we may assume that $\eta\leq 1$,
    \begin{align}\lbeq{|q|-bd}
    |\hat q(k)|\leq|1-\vep(2-\eta)|\vee(1-\vep\eta)\leq1-\vep\eta.
    \end{align}
Therefore, we have
    \begin{align}
    (D^{*2}*q^{*j})(x)&\leq\int_{\hat B_{1/L}^+}\frac{d^dk}{(2\pi)^d}~e^{-cj
    \vep L^2|k|^2}+(1-\vep\eta)^j\int_{\Box_\pi\setminus\hat B_{1/L}^+}\frac
    {d^dk}{(2\pi)^d}~\hat D(k)^2\nn\\
    &=O(\beta)\,(j\vep)^{-d/2}+O(\beta)\,(1-\vep\eta)^j.\lbeq{D*q-bd}
    \end{align}
Substituting \refeq{D*q-bd} into \refeq{S2} and separating the sum of
the first term in \refeq{D*q-bd} into the sum over $N\leq j<\frac{n}2-1$
and the sum over $\frac{n}2-1\leq j\leq n-2$, we obtain
    \begin{align}
    S_2&\leq O(\beta)\,(n\vep-\vep-1)\,(1-\vep)^{\frac{n}2-1}\,\vep\sum_{N\leq
    j<\frac{n}2-1}(j\vep)^{-d/2}\nn\\
    &\qquad+O(\beta)\,(\tfrac{n\vep}2-\vep)^{-d/2}\,\vep^2\sum_{\frac{n}2-1\leq
    j\leq n-2}(n-1-j)\,(1-\vep)^{n-2-j}\nn\\
    &\qquad+O(\beta)\,(1-\vep\eta)^{n-2}\,\vep^2\sum_{N\leq j\leq n-2}(n-1-j)\nn\\
    &\leq O(\beta)\,(n\vep)e^{-n\vep/2}+O(\beta)\,(n\vep)^{-d/2}+O(\beta)\,(n
    \vep)^2\,e^{-\eta n\vep}.\lbeq{S2-bd}
    \end{align}
The proof of \refeq{q-gauss} is completed by
combining \refeq{qbd-small}--\refeq{q-dec2}, \refeq{S1-bd} and \refeq{S2-bd},
and $(n\vep)^2 e^{-n\vep/2}\leq C [1\vee (n\vep)]^{-d/2}$ for all $n\geq 1/\vep$.
\qed

\section{Computation for the spread-out uniform model}
\label{s:uniform}
In this appendix, we compute the model-dependent terms of $\pc-1$ in
\refeq{sawcp-critpt}--\refeq{perc-critpt} when the probability distribution
$D$ is defined as in \refeq{dex}.  Recalling \refeq{U-def}, we have
\begin{align}\lbeq{D*D}
D^{*2}(o)=\frac1{(2L+1)^d-1}=\frac\beta{2^d}+O(\beta L^{-1})=\beta\,
 U^{\star2}(o)+O(\beta L^{-1}).
\end{align}
This relation can be extended as follows:

\begin{prop}\label{prop:DtoU}
Let $D$ be the function defined in \refeq{dex}.  For $\alpha=0,1$,
as $L\to\infty$,
\begin{align}
\sum_{n=3}^\infty(n+1)^\alpha\,D^{*n}(o)&=\beta\sum_{n=3}^\infty
 (n+1)^\alpha\,U^{\star n}(o)+O(\beta L^{-1}),\lbeq{nD*n}\\
\sum_{n=2}^\infty D^{*2n}(o)&=\beta\sum_{n=2}^\infty U^{\star 2n}
 (o)+O(\beta L^{-1}),\lbeq{D*2n}
\end{align}
where $d>4+2\alpha$ in \refeq{nD*n} and $d>4$ in \refeq{D*2n}.
\end{prop}

Theorem \ref{thm:uniform} is an immediate consequence of \refeq{D*D} and
Proposition~\ref{prop:DtoU}.  Note further that the coefficients of
$\beta$ in \refeq{nD*n}--\refeq{D*2n} are bounded if $d>2+2\alpha$ and
$d>2$, respectively, which suggests that also for $d=3+2\alpha$ and
$d=4+2\alpha$ the leading order contributions should be given by the
first terms in \refeq{nD*n}--\refeq{D*2n}.

\begin{proof}
For $x\in\Zd$, define
\begin{align}
D_o(x)=\frac{\indic{\|x\|_\infty\leq L}}{(2L+1)^d},
\end{align}
to be a regularized version of $D$ in \refeq{dex}.  It is obvious that
$D_o^{*2}(o)$ satisfies the same estimate as in \refeq{D*D}.  We note that
\begin{align}
D^{*m}(o)-D_o^{*m}(o)&=\sum_{j=1}^m\big((D-D_o)*D^{*(j-1)}*D_o^{*(m-j)}
 \big)(o)\nn\\
&=\sum_{j=1}^m\Bigg[\sum_{x:0<\|x\|_\infty\leq L}\frac{\big(D^{*(j-1)}
 *D_o^{*(m-j)}\big)(x)} {(2L+1)^d[(2L+1)^d-1]}-\frac{\big(D^{*(j-1)}*
 D_o^{*(m-j)}\big)(o)}{(2L+1)^d}\Bigg]\nn\\
&=\frac1{(2L+1)^d}\sum_{j=1}^m\Big[\big(D^{*j}*D_o^{*(m-j)}\big)(o)-
 \big(D^{*(j-1)}*D_o^{*(m-j)}\big)(o)\Big].\lbeq{D-Do}
\end{align}
By this identity and the fact that $D_o^{*n}(x)$ also satisfies
\refeq{d-gauss}, we can approximate the expressions in the left-hand side
of \refeq{nD*n}--\refeq{D*2n} by the corresponding expressions defined in
terms of $D_o$ instead of $D$, up to $O(\beta^2)$ when $d>4+2\alpha$ and
$d>4$, respectively.  For example, for \refeq{nD*n} with $\alpha=0$, we
use \refeq{d-gauss} to obtain
\begin{align}
&\bigg|\sum_{n=3}^\infty D^{*n}(o)-\sum_{n=3}^\infty D_o^{*n}(o)\bigg|\nn\\
&\qquad\leq\frac1{(2L+1)^d}\sum_{j=1}^\infty\;\sum_{n=j\vee3}^\infty\Big[
 \big(D^{*j}*D_o^{*(n-j)}\big)(o)+\big(D^{*(j-1)}*D_o^{*(n-j)}\big)(o)
 \Big]\nn\\
&\qquad=\frac1{(2L+1)^d}\Big((D+\delta_o)*\big(D_o^{*2}+D*D_o+D^{*2}*G\big)
 *G_o\Big)(o)=O(\beta^2),\lbeq{sumD-sumDo}
\end{align}
where $\delta_o(x)=\delta_{o,x}$ and $G_o(x)=\sum_{n=0}^\infty D_o^{*n}(x)$.
Therefore, to prove Proposition~\ref{prop:DtoU}, it suffices to show that,
for $\alpha=0,1$,
\begin{align}
\sum_{n=3}^\infty(n+1)^\alpha\,D_o^{*n}(o)&=\beta\sum_{n=3}^\infty
 (n+1)^\alpha\,U^{\star n}(o)+O(\beta L^{-1}),\lbeq{nDo*n}\\
\sum_{n=2}^\infty D_o^{*2n}(o)&=\beta\sum_{n=2}^\infty U^{\star 2n}
 (o)+O(\beta L^{-1}).\lbeq{Do*2n}
\end{align}

We prove \refeq{nDo*n} for $\alpha=0$ by comparing the Fourier transform
of the first term in the right-hand side of \refeq{nDo*n}, i.e.,
\begin{align}\lbeq{Usum-fourier}
\beta\sum_{n=3}^\infty U^{\star n}(o)=\beta\int_{\mR^d}\frac{d^dk}
 {(2\pi)^d}~\frac{\hat U(k)^3}{1-\hat U(k)},
\end{align}
with the Fourier transform of the left-hand side of \refeq{nDo*n}, i.e.,
\begin{align}
\sum_{n=3}^\infty D_o^{*n}(o)=\int_{\Box_\pi}\frac{d^dk}{(2\pi)^d}~\frac
 {\hat D_o(k)^3}{1-\hat D_o(k)}&=\beta_o\int_{\Box_{
 (L+\frac12)\pi}}\frac{d^dk}{(2\pi)^d}~\frac{\hat D_o\big(\frac{k}{L+\frac
 12}\big)^3}{1-\hat D_o\big(\frac{k}{L+\frac12}\big)},\lbeq{Dosum-fourier}
\end{align}
where $\Box_\ell=[-\ell,\ell]^d$ and
\begin{align}
\beta_o=(L+\tfrac12)^{-d}=\beta+O(\beta L^{-1}),\lbeq{betao}
\end{align}
and also
\begin{align}
\hat U(k)&=\int_{\mR^d}d^dx~U(x)~e^{ik\cdot x}=\prod_{j=1}^d\frac{\sin
 k_j}{k_j},\lbeq{U-fourier}\\
\hat D_o(k)&=\sum_{x\in\Zd}D_o(x)~e^{ik\cdot x}=\prod_{j=1}^d\frac{\sin
 [(2L+1)\frac{k_j}2]}{(2L+1)\sin\frac{k_j}2}=\frac{\hat U((L+\tfrac12)k)}
 {\hat U(\frac{k}2)}.\lbeq{Do-fourier}
\end{align}
The simple product formula \refeq{Do-fourier} is the main advantage of
using $D_o$ instead of $D$.  It follows from \refeq{Do-fourier} that
$\hat D_o\big(\frac{k}{L+\frac12}\big)=\hat U(k)/\hat U(\frac{k}{2L+1})$,
which approximates $\hat U(k)$ for large $L$.  We write
\begin{align}\lbeq{I12}
\beta\sum_{n=3}^\infty U^{\star n}(o)-\sum_{n=3}^\infty D_o^{*n}(o)=
 (\beta-\beta_o)\sum_{n=3}^\infty U^{\star n}(o)+\beta_o(I_1+I_2),
\end{align}
where, by \refeq{betao}, the first term is $O(\beta L^{-1})$, and
\begin{align}\lbeq{I1I2}
I_1=\int_{\mR^d\setminus\Box_{(L+\frac12)\pi}}\frac{d^dk}{(2\pi)^d}~
 \frac{\hat U(k)^3}{1-\hat U(k)},&&
I_2=\int_{\Box_{(L+\frac12)\pi}}\frac{d^dk}{(2\pi)^d}~\Bigg[\frac{
 \hat U(k)^3}{1-\hat U(k)}-\frac{\hat D_o\big(\frac{k}{L+\frac12}
 \big)^3}{1-\hat D_o\big(\frac{k}{L+\frac12}\big)}\Bigg].
\end{align}
We prove below that each $I_j$ is $O(L^{-1})$ if $d>2$.  This suffices
to prove \refeq{nDo*n} for $\alpha=0$.  If fact, we will prove that
each $I_j$ is $O(L^{-2}\log L)$ if $d>2$, which also identifies the
coefficient of $\beta L^{-1}$.

To estimate each $I_j$, we use the following properties
of $\hat U(k)$ and $\hat D_o\big(\frac{k}{L+\frac12}\big)$ that follow
from the standard estimates for the trigonometric functions: for any $k$,
\begin{align}\lbeq{U-summary}
|\hat U(k)|\leq\prod_{j=1}^d(1\vee|k_j|)^{-1},&&
1-\hat U(k)&\geq c_1(1\wedge|k|^2),
\end{align}
and for $k\in\Box_{(L+\frac12)\pi}$,
\begin{gather}
\big|\hat D_o\big(\tfrac{k}{L+\frac12}\big)\big|\leq c_2\prod_{j=1}^d
 (1\vee|k_j|)^{-1},\qquad\qquad
1-\hat D_o\big(\tfrac{k}{L+\frac12}\big)\geq c_3(1\wedge|k|^2),
 \lbeq{Do-summary}\\[5pt]
\big|\hat U(k)-\hat D_o\big(\tfrac{k}{L+\frac12}\big)\big|\leq c_4
 L^{-2}|\hat U(k)|\,|k|^2,\lbeq{U-Do-summary}
\end{gather}
where each $c_i\in(0,\infty)$ is independent of $L$ and $k$.  To see
the first inequality in \refeq{U-summary}, we only need to use
$|\frac{\sin r}r|\leq(1\vee r)^{-1}$ for any $r$.  For the first
inequality in \refeq{Do-summary}, we recall $\hat D_o\big(\frac{k}
{L+\frac12}\big)=\hat U(k)/\hat U(\frac{k}{2L+1})$ and use
$\frac{\sin r}r\geq\frac2\pi$ for any $r\in[-\frac\pi2,\frac\pi2]$,
so that $|\hat U(\frac{k}{2L+1})^{-1}|\leq(\frac\pi2)^d$, and then
use the bound on $|\hat U(k)|$ in \refeq{U-summary}.  The second
inequalities in \refeq{U-summary}--\refeq{Do-summary} follow from
\refeq{U-fourier} and \refeq{Dprop}, respectively.  Finally, for
\refeq{U-Do-summary}, we again use $\frac2\pi\leq\frac{\sin r}r\leq1$
for any $r\in[-\frac\pi2,\frac\pi2]$ to obtain
\begin{align}
\big|\hat U(k)-\hat D_o\big(\tfrac{k}{L+\frac12}\big)\big|&\leq\Big(\frac
 \pi2\Big)^{\!d}\,|\hat U(k)|\,\bigg|\prod_{j=1}^d\frac{\sin\frac{k_j}{2L
 +1}}{\frac{k_j}{2L+1}}-1\bigg|\leq\Big(\frac\pi2\Big)^{\!d}\,|\hat U(k)|
 \sum_{i=1}^d\bigg|\frac{\sin\frac{k_i}{2L+1}}{\frac{k_i}{2L+1}}-1\bigg|,
\end{align}
which is bounded by $|\hat U(k)|\,O(|k|^2L^{-2})$, using
$0\leq\frac{\sin r}r-1+\frac{r^2}{3!}\leq\frac{r^2}{5!}(\frac\pi2)^2$
for any $r\in[-\frac\pi2,\frac\pi2]$.  This completes the proof of
\refeq{U-summary}--\refeq{U-Do-summary}.

First, we consider $I_1$.  Using \refeq{U-summary}, we obtain
\begin{align}
|I_1|\leq c\int_{\mR^d\setminus\Box_{(L+\frac12)\pi}}d^dk~|\hat U(k)|^3
 \leq c'\sum_{i=1}^d\bigg[\int_{(L+\frac12)\pi}^\infty\frac{dk_i}{k_i^3}
 \bigg]\bigg[\prod_{j\ne i}\int_{-\infty}^\infty\frac{dk_j}{(1\vee|k_j|)^3}
 \bigg]=O(L^{-2}).\lbeq{I1-bd}
\end{align}
For $I_2$, we use
\begin{align}
\bigg|\frac{u^3}{1-u}-\frac{d^3}{1-d}\bigg|=\frac{|u-d|\,|u^2+ud+d^2
 -ud(u+d)|}{(1-u)(1-d)}\leq\frac{|u-d|\,(u^2+3|ud|+d^2)}{(1-u)(1-d)},
\end{align}
with $u=\hat U(k)$ and $d=\hat D_o\big(\frac{k}{L+\frac12}\big)$.
Using \refeq{U-summary}--\refeq{U-Do-summary}, we obtain
\begin{align}
|I_2|\leq cL^{-2}\int_{\Box_{(L+\frac12)\pi}}d^dk~\frac{|k|^2\prod_{j=
 1}^d(1\vee|k_j|)^{-3}}{(1\wedge|k|^2)^2}=O(L^{-2}\log L),\lbeq{I2-bd}
\end{align}
for $d>2$.
The proof of \refeq{nDo*n} for $\alpha=0$ is completed by \refeq{I12},
\refeq{I1-bd} and \refeq{I2-bd}.

The same strategy explained above also applies to the proof of \refeq{nD*n}
for $\alpha=1$ and \refeq{D*2n}, using the following expressions:
\begin{align}
\sum_{n=3}^\infty(n+1)\,D_o^{*n}(o)&=\beta_o\int_{\Box_{(L+\frac12)\pi}}
 \frac{d^dk}{(2\pi)^d}~\frac{\hat D_o\big(\frac{k}{L+\frac12}\big)^3
 \big[4-3\hat D_o\big(\frac{k}{L+\frac12}\big)\big]}{\big[1-\hat D_o
 \big(\frac{k}{L+\frac12}\big)\big]^2},\\
\sum_{n=2}^\infty D_o^{*2n}(o)&=\beta_o\int_{\Box_{(L+\frac12)\pi}}
 \frac{d^dk}{(2\pi)^d}~\frac{\hat D_o\big(\frac{k}{L+\frac12}\big)^4}
 {1-\hat D_o\big(\frac{k}{L+\frac12}\big)^2}.
\end{align}
This completes the proof of Proposition~\ref{prop:DtoU}.
\end{proof}

\section*{Acknowledgements}
The work of AS was supported in part by NSERC of Canada.  The work of RvdH
and AS was supported in part by Netherlands Organisation for Scientific
Research (NWO). 
This project was initiated during an extensive visit of RvdH to the
University of British Columbia, Vancouver, Canada. We thank Gordon~Slade
for comments on a preliminary version of the paper, and for pointing us
to the appropriate bound in \cite{HS90a} that implies Lemma \ref{lem-T}.


\begin{thebibliography}{99}

\bibitem{an84}M. Aizenman and C. M. Newman.
\newblock Tree graph inequalities and critical behavior in percolation models.
\newblock {\it J. Statist. Phys.} {\bf 36} (1984): 107--143.

\bibitem{ba91}D. J. Barsky and M. Aizenman.
\newblock Percolation critical exponents under the triangle condition.
\newblock {\it Ann. Probab.} {\bf 19} (1991): 1520--1536.



\bibitem{bg91}C. Bezuidenhout and G. Grimmett.
\newblock Exponential decay for subcritical contact and percolation processes.
\newblock {\it Ann. Probab.} {\bf 19} (1991): 984--1009.

\bibitem{bds89}M. Bramson, R. Durrett and G. Swindle.
\newblock Statistical mechanics of crabgrass.
\newblock {\it Ann. Probab.} {\bf 17} (1989): 444--481.


\bibitem{BCHSS04b}
C.~Borgs, J.T. Chayes, R.~van~der Hofstad, G.~Slade, and J.~Spencer.
\newblock Random subgraphs of finite graphs: {II}. {The} lace expansion
and the triangle condition.
\newblock Preprint, (2003).

\bibitem{cd83}J. T. Cox and R. Durrett.
\newblock Oriented percolation in dimensions $d\geq4$: bounds and asymptotic
formulars.
\newblock {\it Math. Proc. Cambridge Philos. Soc.} {\bf 93} (1983): 151-162.



\bibitem{dp99}R. Durrett and E. Perkins.
\newblock Rescaled contact processes converge to super-Brownian motion in two
or more dimensions.
\newblock {\it Probab. Th. Rel. Fields} {\bf 114} (1999): 309--399.







\bibitem{g99}G. Grimmett.
\newblock {\it Percolation}.  Springer, Berlin (1999).

\bibitem{gh01}G. Grimmett and P. Hiemer.
\newblock Directed percolation and random walk.
\newblock {\it In and Out of Equilibrium} (ed. V. Sidoravicius).
Birkh\"auser (2002): 273-297.

\bibitem{HHS01a}
T.~Hara, R.~van~der Hofstad, and G.~Slade.
\newblock Critical two-point functions and the lace expansion for spread-out
  high-dimensional percolation and related models.
\newblock {\em Ann.\ Probab.}, {\bf 31} (2003): 349-408.

\bibitem{HS90a}
T.~Hara and G.~Slade.
\newblock Mean-field critical behaviour for percolation in high dimensions.
\newblock {\em Commun. Math. Phys.}, {\bf 128} (1990): 333--391.


\bibitem{hs95}
T.~Hara and G.~Slade.
\newblock The self-avoiding-walk and percolation critical points in high dimensions.
\newblock {\it Combin. Probab. Comput.}, {\bf 4} (1995): 197--215.

%

\bibitem{h03}R. van der Hofstad.
\newblock The derivative of the lace expansion coefficients for
unoriented percolation.
\newblock Unpublished document (2003).

%

\bibitem{hsa04}R. van der Hofstad and A. Sakai.
\newblock Gaussian scaling for the critical spread-out
contact process above the upper critical dimension.
\newblock {\it Preprint} (2003). To appear in {\it Electr.\ Journ.\ Probab.}

\bibitem{hs02}R. van der Hofstad and G. Slade.
\newblock A generalised inductive approach to the lace expansion.
\newblock {\it Probab. Th. Rel. Fields} {\bf 122} (2002): 389--430.

\bibitem{hs01}R. van der Hofstad and G. Slade.
\newblock Convergence of critical oriented percolation to super-Brownian motion
above 4+1 dimensions.
\newblock {\it Ann. Inst. H. Poincar\'e Probab. Statist.} {\bf 39} (2003): 413--485.

\bibitem{hs03}R. van der Hofstad and G. Slade.
\newblock The lace expansion on a tree with application to networks of
self-avoiding walks.
\newblock {\it Adv. Appl. Math.} {\bf 30} (2003): 471--528.

\bibitem{HS04a}R. van der Hofstad and G. Slade.
\newblock Expansion in $n^{-1}$ for percolation critical values on the
$n$-cube and ${\mathbb Z}^n$: the first three terms.
\newblock {\it Preprint} (2004). To appear in {\it Combin. Probab. Comput.}

\bibitem{HS04b}R. van der Hofstad and G. Slade.
\newblock Asymptotic expansions in $n^{-1}$ for percolation critical values
on the $n$-cube and ${\mathbb Z}^n$.
\newblock {\it Preprint} (2004).

\bibitem{k80}H. Kesten.
\newblock The critical probability of bond percolation on the square
lattice equals $\frac12$.
\newblock {\it Commun. Math. Phys.} {\bf 74} (1980): 41--59.

\bibitem{Ligg99}T. Liggett.
\newblock {\em Stochastic Interacting Systems: Contact, Voter and Exclusion Processes}.
\newblock Springer, Berlin (1999).

\bibitem{l97}T. Liggett.
\newblock Stochastic models of interacting systems.
\newblock {\it Ann. Probab.} {\bf 25} (1997): 1--29.


\bibitem{ms93}N. Madras and G. Slade.
\newblock {\it The Self-Avoiding Walk}.
\newblock Birkh\"auser, Boston (1993).

\bibitem{ny93}B. G. Nguyen and W.-S. Yang.
\newblock Triangle condition for oriented percolation in high dimensions.
\newblock {\it Ann. Probab.} {\bf 21} (1993): 1809--1844.


\bibitem{p94}M. D. Penrose.
\newblock Self-avoiding walks and trees in spread-out lattices.
\newblock {\it J. Statist. Phys.} {\bf 77} (1994): 3--15.



\bibitem{s01}A. Sakai.
\newblock Mean-field critical behavior for the contact process.
\newblock {\it J. Statist. Phys.} {\bf 104} (2001): 111--143.






\end{thebibliography}
\end{document}